# THE FORMAL DEFINITION OF REFERENCE PRIORS

By James O. Berger,[1] José M. Bernardo[2] and Dongchu Sun[3]

*Duke University, Universitat de València and University of Missouri-Columbia*

Reference analysis produces objective Bayesian inference, in the sense that inferential statements depend only on the assumed model and the available data, and the prior distribution used to make an inference is least informative in a certain information-theoretic sense. Reference priors have been rigorously defined in specific contexts and heuristically defined in general, but a rigorous general definition has been lacking. We produce a rigorous general definition here and then show how an explicit expression for the reference prior can be obtained under very weak regularity conditions. The explicit expression can be used to derive new reference priors both analytically and numerically.

## 1. Introduction and notation.

1.1. *Background and goals.* There is a considerable body of conceptual and theoretical literature devoted to identifying appropriate procedures for the formulation of objective priors; for relevant pointers see Section 5.6 in Bernardo and Smith [13], Datta and Mukerjee [20], Bernardo [11], Berger [3], Ghosh, Delampady and Samanta [23] and references therein. Reference analysis, introduced by Bernardo [10] and further developed by Berger and Bernardo [4, 5, 6, 7], and Sun and Berger [42], has been one of the most utilized approaches to developing objective priors; see the references in Bernardo [11].

Reference analysis uses information-theoretical concepts to make precise the idea of an objective prior which should be maximally dominated by the

Received March 2007; revised December 2007.
[1]Supported by NSF Grant DMS-01-03265.
[2]Supported by Grant MTM2006-07801.
[3]Supported by NSF Grants SES-0351523 and SES-0720229.

data, in the sense of maximizing the missing information (to be precisely defined later) about the parameter. The original formulation of reference priors in the paper by Bernardo [10] was largely informal. In continuous one parameter problems, heuristic arguments were given to justify an explicit expression in terms of the expectation under sampling of the logarithm of the asymptotic posterior density, which reduced to Jeffreys prior (Jeffreys [31, 32]) under asymptotic posterior normality. In multiparameter problems it was argued that one should not maximize the joint missing information but proceed sequentially, thus avoiding known problems such as marginalization paradoxes. Berger and Bernardo [7] gave more precise definitions of this sequential reference process, but restricted consideration to continuous multiparameter problems under asymptotic posterior normality. Clarke and Barron [17] established regularity conditions under which joint maximization of the missing information leads to Jeffreys multivariate priors. Ghosal and Samanta [27] and Ghosal [26] provided explicit results for reference priors in some types of nonregular models.

This paper has three goals.

GOAL 1. *Make precise the definition of the reference prior.* This has two different aspects.

- Applying Bayes theorem to improper priors is not obviously justifiable. Formalizing when this is legitimate is desirable, and is considered in Section 2.
- Previous attempts at a general definition of reference priors have had heuristic features, especially in situations in which the reference prior is improper. Replacing the heuristics with a formal definition is desirable, and is done in Section 3.

GOAL 2. *Present a simple constructive formula for a reference prior.* Indeed, for a model described by density $p(\mathbf{x} \mid \theta)$, where $\mathbf{x}$ is the complete data vector and $\theta$ is a continuous unknown parameter, the formula for the reference prior, $\pi(\theta)$, will be shown to be

$$\pi(\theta) = \lim_{k \to \infty} \frac{f_k(\theta)}{f_k(\theta_0)},$$

$$f_k(\theta) = \exp\left\{ \int p(\mathbf{x}^{(k)} \mid \theta) \log[\pi^*(\theta \mid \mathbf{x}^{(k)})] \, d\mathbf{x}^{(k)} \right\},$$

where $\theta_0$ is an interior point of the parameter space $\Theta$, $\mathbf{x}^{(k)} = \{\mathbf{x}_1, \ldots, \mathbf{x}_k\}$ stands for $k$ conditionally independent replications of $\mathbf{x}$, and $\pi^*(\theta \mid \mathbf{x}^{(k)})$ is the posterior distribution corresponding to some fixed, largely arbitrary prior $\pi^*(\theta)$.

The interesting thing about this expression is that it holds (under mild conditions) for any type of continuous parameter model, regardless of the asymptotic nature of the posterior. This formula is established in Section 4.1, and various illustrations of its use are given.

A second use of the expression is that it allows straightforward computation of the reference prior numerically. This is illustrated in Section 4.2 for a difficult nonregular problem and for a problem for which analytical determination of the reference prior seems very difficult.

GOAL 3. *To make precise the most common practical rationale for use of improper objective priors, which proceeds as follows:*

- In reality, we are always dealing with bounded parameters so that the real parameter space should, say, be some compact set $\Theta_0$.
- It is often only known that the bounds are quite large, in which case it is difficult to accurately ascertain which $\Theta_0$ to use.
- This difficulty can be surmounted if we can pass to the unbounded space $\Theta$ and show that the analysis on this space would yield essentially the same answer as the analysis on any very large compact $\Theta_0$.

Establishing that the analysis on $\Theta$ is a good approximation from the reference theory viewpoint requires establishing two facts:

1. The reference prior distribution on $\Theta$, when restricted to $\Theta_0$, is the reference prior on $\Theta_0$.
2. The reference posterior distribution on $\Theta$ is an appropriate limit of the reference posterior distributions on an increasing sequence of compact sets $\{\Theta_i\}_{i=1}^{\infty}$ converging to $\Theta$.

Indicating how these two facts can be verified is the third goal of the paper.

1.2. *Notation.* Attention here is limited mostly to one parameter problems with a continuous parameter, but the ideas are extendable to the multiparameter case through the sequential scheme of Berger and Bernardo [7].

It is assumed that probability distributions may be described through probability density functions, either in respect to Lebesgue measure or counting measure. No distinction is made between a random quantity and the particular values that it may take. Bold italic roman fonts are used for observable random vectors (typically data) and italic greek fonts for unobservable random quantities (typically parameters); lower case is used for variables and upper case calligraphic for their domain sets. Moreover, the standard mathematical convention of referring to functions, say $f_{\mathbf{x}}$ and $g_{\mathbf{x}}$ of $\mathbf{x} \in \mathcal{X}$, respectively by $f(\mathbf{x})$ and $g(\mathbf{x})$, will be used throughout. Thus, the conditional probability density of data $\mathbf{x} \in \mathcal{X}$ given $\theta$ will be represented by $p(\mathbf{x} \mid \theta)$, with $p(\mathbf{x} \mid \theta) \geq 0$ and $\int_{\mathcal{X}} p(\mathbf{x} \mid \theta) \, d\mathbf{x} = 1$, and the reference posterior distribution of $\theta \in \Theta$ given $\mathbf{x}$ will be represented by $\pi(\theta \mid \mathbf{x})$, with $\pi(\theta \mid \mathbf{x}) \geq 0$ and $\int_\Theta \pi(\theta \mid \mathbf{x}) \, d\theta = 1$. This admittedly imprecise notation will greatly simplify the exposition. If the random vectors are discrete, these functions naturally become probability mass functions, and integrals over their values become sums. Density functions of specific distributions are denoted by appropriate names. Thus, if $x$ is an observable random quantity with a normal distribution of mean $\mu$ and variance $\sigma^2$, its probability density function will be denoted $\mathrm{N}(x \mid \mu, \sigma^2)$; if the posterior distribution of $\lambda$ is Gamma with mean $a/b$ and variance $a/b^2$, its probability density function will be denoted $\mathrm{Ga}(\lambda \mid a, b)$. The indicator function on a set $C$ will be denoted by $\mathbf{1}_C$.

Reference prior theory is based on the use of logarithmic divergence, often called the *Kullback–Leibler divergence*.

DEFINITION 1. The logarithmic divergence of a probability density $\tilde{p}(\theta)$ of the random vector $\theta \in \Theta$ from its true probability density $p(\theta)$, denoted by $\kappa\{\tilde{p} \mid p\}$, is

$$\kappa\{\tilde{p} \mid p\} = \int_\Theta p(\theta) \log \frac{p(\theta)}{\tilde{p}(\theta)} \, d\theta,$$

provided the integral (or the sum) is finite.

The properties of $\kappa\{\tilde{p} \mid p\}$ have been extensively studied; pioneering works include Gibbs [22], Shannon [38], Good [24, 25], Kullback and Leibler [35], Chernoff [15], Jaynes [29, 30], Kullback [34] and Csiszar [18, 19].

DEFINITION 2 (Logarithmic convergence). A sequence of probability density functions $\{p_i\}_{i=1}^\infty$ converges logarithmically to a probability density $p$ if, and only if, $\lim_{i \to \infty} \kappa(p \mid p_i) = 0$.

## 2. Improper and permissible priors.

2.1. *Justifying posteriors from improper priors.* Consider a model $\mathcal{M} = \{p(\mathbf{x} \mid \theta), \mathbf{x} \in \mathcal{X}, \theta \in \Theta\}$ and a strictly positive prior function $\pi(\theta)$. (We restrict attention to strictly positive functions because any believably objective prior would need to have strictly positive density, and this restriction eliminates many technical details.) When $\pi(\theta)$ is improper, so that $\int_\Theta \pi(\theta) \, d\theta$ diverges, Bayes theorem no longer applies, and the use of the formal posterior density

$$\pi(\theta \mid \mathbf{x}) = \frac{p(\mathbf{x} \mid \theta) \pi(\theta)}{\int_\Theta p(\mathbf{x} \mid \theta) \pi(\theta) \, d\theta} \tag{2.1}$$

must be justified, even when $\int_\Theta p(\mathbf{x} \mid \theta)\pi(\theta) \, d\theta < \infty$ so that $\pi(\theta \mid \mathbf{x})$ is a proper density.

The most convincing justifications revolve around showing that $\pi(\theta \mid \mathbf{x})$ is a suitable limit of posteriors obtained from proper priors. A variety of versions of such arguments exist; cf. Stone [40, 41] and Heath and Sudderth [28]. Here, we consider approximations based on restricting the prior to an increasing sequence of compact sets and using logarithmic convergence to define the limiting process. The main motivation is, as mentioned in the introduction, that objective priors are often viewed as being priors that will yield a good approximation to the analysis on the "true but difficult to specify" large bounded parameter space.

DEFINITION 3 (Approximating compact sequence). Consider a parametric model $\mathcal{M} = \{p(\mathbf{x} \mid \theta), \mathbf{x} \in \mathcal{X}, \theta \in \Theta\}$ and a strictly positive continuous function $\pi(\theta)$, $\theta \in \Theta$, such that, for all $\mathbf{x} \in \mathcal{X}$, $\int_\Theta p(\mathbf{x} \mid \theta)\pi(\theta) \, d\theta < \infty$. An approximating compact sequence of parameter spaces is an increasing sequence of compact subsets of $\Theta$, $\{\Theta_i\}_{i=1}^\infty$, converging to $\Theta$. The corresponding sequence of posteriors with support on $\Theta_i$, defined as $\{\pi_i(\theta \mid \mathbf{x})\}_{i=1}^\infty$, with $\pi_i(\theta \mid \mathbf{x}) \propto p(\mathbf{x} \mid \theta)\pi_i(\theta)$, $\pi_i(\theta) = c_i^{-1}\pi(\theta)\mathbf{1}_{\Theta_i}$ and $c_i = \int_{\Theta_i} \pi(\theta) \, d\theta$, is called the approximating sequence of posteriors to the formal posterior $\pi(\theta \mid \mathbf{x})$.

Notice that the renormalized restrictions $\pi_i(\theta)$ of $\pi(\theta)$ to the $\Theta_i$ are proper [because the $\Theta_i$ are compact and $\pi(\theta)$ is continuous]. The following theorem shows that the posteriors resulting from these proper priors do converge, in the sense of logarithmic convergence, to the posterior $\pi(\theta \mid \mathbf{x})$.

THEOREM 1. *Consider model $\mathcal{M} = \{p(\mathbf{x} \mid \theta), \mathbf{x} \in \mathcal{X}, \theta \in \Theta\}$ and a strictly positive continuous function $\pi(\theta)$, such that $\int_\Theta p(\mathbf{x} \mid \theta)\pi(\theta) \, d\theta < \infty$, for all $\mathbf{x} \in \mathcal{X}$. For any approximating compact sequence of parameter spaces, the corresponding approximating sequence of posteriors converges logarithmically to the formal posterior $\pi(\theta \mid \mathbf{x}) \propto p(\mathbf{x} \mid \theta)\pi(\theta)$.*

PROOF. To prove that $\kappa\{\pi(\cdot \mid \mathbf{x}) \mid \pi_i(\cdot \mid \mathbf{x})\}$ converges to zero, define the predictive densities $p_i(\mathbf{x}) = \int_{\Theta_i} p(\mathbf{x} \mid \theta)\pi_i(\theta) \, d\theta$ and $p(\mathbf{x}) = \int_\Theta p(\mathbf{x} \mid \theta)\pi(\theta) \, d\theta$ (which has been assumed to be finite). Using for the posteriors the expressions provided by Bayes theorem yields

$$\begin{aligned}
\int_{\Theta_i} \pi_i(\theta \mid \mathbf{x}) \log \frac{\pi_i(\theta \mid \mathbf{x})}{\pi(\theta \mid \mathbf{x})} \, d\theta &= \int_{\Theta_i} \pi_i(\theta \mid \mathbf{x}) \log \frac{p(\mathbf{x})\pi_i(\theta)}{p_i(\mathbf{x})\pi(\theta)} \, d\theta \\
&= \int_{\Theta_i} \pi_i(\theta \mid \mathbf{x}) \log \frac{p(\mathbf{x})}{p_i(\mathbf{x})c_i} \, d\theta \\
&= \log \frac{p(\mathbf{x})}{p_i(\mathbf{x})c_i} = \log \frac{\int_\Theta p(\mathbf{x} \mid \theta)\pi(\theta) \, d\theta}{\int_{\Theta_i} p(\mathbf{x} \mid \theta)\pi(\theta) \, d\theta}.
\end{aligned}$$

But the last expression converges to zero if, and only if,

$$\lim_{i \to \infty} \int_{\Theta_i} p(\mathbf{x} \mid \theta) \pi(\theta) \, d\theta = \int_{\Theta} p(\mathbf{x} \mid \theta) \pi(\theta) \, d\theta,$$

and this follows from the monotone convergence theorem. □

It is well known that logarithmic convergence implies convergence in $L_1$ which implies uniform convergence of probabilities, so Theorem 1 could, at first sight, be invoked to justify the formal use of virtually any improper prior in Bayes theorem. As illustrated below, however, logarithmic convergence of the approximating posteriors is not necessarily good enough.

EXAMPLE 1 (Fraser, Monette and Ng [21]). Consider the model, with both discrete data and parameter space,

$$\mathcal{M} = \{p(x \mid \theta) = 1/3, \ x \in \{[\theta/2], 2\theta, 2\theta + 1\}, \ \theta \in \{1, 2, \ldots\}\},$$

where $[u]$ denotes the integer part of $u$, and $[1/2]$ is separately defined as 1. Fraser, Monnete and Ng [21] show that the naive improper prior $\pi(\theta) = 1$ produces a posterior $\pi(\theta \mid x) \propto p(x \mid \theta)$ which is strongly inconsistent, leading to credible sets for $\theta$ given by $\{2x, 2x + 1\}$ which have posterior probability $2/3$ but frequentist coverage of only $1/3$ for all $\theta$ values. Yet, choosing the natural approximating sequence of compact sets $\Theta_i = \{1, \ldots, i\}$, it follows from Theorem 1 that the corresponding sequence of posteriors converges logarithmically to $\pi(\theta \mid x)$.

The difficulty shown by Example 1 lies in the fact that logarithmic convergence is only pointwise convergence for given $\mathbf{x}$, which does not guarantee that the approximating posteriors are accurate in any global sense over $\mathbf{x}$. For that we turn to a stronger notion of convergence.

DEFINITION 4 (Expected logarithmic convergence of posteriors). Consider a parametric model $\mathcal{M} = \{p(\mathbf{x} \mid \theta), \mathbf{x} \in \mathcal{X}, \theta \in \Theta\}$, a strictly positive continuous function $\pi(\theta)$, $\theta \in \Theta$ and an approximating compact sequence $\{\Theta_i\}$ of parameter spaces. The corresponding sequence of posteriors $\{\pi_i(\theta \mid \mathbf{x})\}_{i=1}^{\infty}$ is said to be expected logarithmically convergent to the formal posterior $\pi(\theta \mid \mathbf{x})$ if

$$(2.2) \qquad \lim_{i \to \infty} \int_{\mathcal{X}} \kappa\{\pi(\cdot \mid \mathbf{x}) \mid \pi_i(\cdot \mid \mathbf{x})\} p_i(\mathbf{x}) \, d\mathbf{x} = 0,$$

where $p_i(\mathbf{x}) = \int_{\Theta_i} p(\mathbf{x} \mid \theta) \pi_i(\theta) \, d\theta$.

This notion was first discussed (in the context of reference priors) in Berger and Bernardo [7], and achieves one of our original goals: A prior distribution satisfying this condition will yield a posterior that, on average over **x**, is a good approximation to the proper posterior that would result from restriction to a large compact subset of the parameter space.

To some Bayesians, it might seem odd to worry about averaging the logarithmic discrepancy over the sample space but, as will be seen, reference priors are designed to be "noninformative" for a specified model, the notion being that repeated use of the prior with that model will be successful in practice.

EXAMPLE 2 (Fraser, Monette and Ng [21] continued). In Example 1, the discrepancies $\kappa\{\pi(\cdot \mid \mathbf{x}) \mid \pi_i(\cdot \mid \mathbf{x})\}$ between $\pi(\theta \mid \mathbf{x})$ and the posteriors derived from the sequence of proper priors $\{\pi_i(\theta)\}_{i=1}^\infty$ converged to zero. However, Berger and Bernardo [7] shows that $\int_{\mathcal{X}} \kappa\{\pi(\cdot \mid \mathbf{x}) \mid \pi_i(\cdot \mid \mathbf{x})\} p_i(\mathbf{x}) \, d\mathbf{x} \to \log 3$ as $i \to \infty$, so that the expected logarithmic discrepancy does not go to zero. Thus, the sequence of proper priors $\{\pi_i(\theta) = 1/i, \theta \in \{1, \ldots, i\}\}_{i=1}^\infty$ does not provide a good global approximation to the formal prior $\pi(\theta) = 1$, providing one explanation of the paradox found by Fraser, Monette and Ng [21].

Interestingly, for the improper prior $\pi(\theta) = 1/\theta$, the approximating compact sequence considered above can be shown to yield posterior distributions that expected logarithmically converge to $\pi(\theta \mid x) \propto \theta^{-1} p(x \mid \theta)$, so that this is a good candidate objective prior for the problem. It is also shown in Berger and Bernardo [7] that this prior has posterior confidence intervals with the correct frequentist coverage.

Two potential generalizations are of interest. Definition 4 requires convergence only with respect to one approximating compact sequence of parameter spaces. It is natural to wonder what happens for other such approximating sequences. We suspect, but have been unable to prove in general, that convergence with respect to one sequence will guarantee convergence with respect to any sequence. If true, this makes expected logarithmic convergence an even more compelling property.

Related to this is the possibility of allowing not just an approximating series of priors based on truncation to compact parameter spaces, but instead allowing any approximating sequence of priors. Among the difficulties in dealing with this is the need for a better notion of divergence that is symmetric in its arguments. One possibility is the symmetrized form of the logarithmic divergence in Bernardo and Rueda [12], but the analysis is considerably more difficult.

2.2. *Permissible priors.* Based on the previous considerations, we restrict consideration of possibly objective priors to those that satisfy the expected logarithmic convergence condition, and formally define them as follows. (Recall that **x** represents the entire data vector.)

DEFINITION 5. A strictly positive continuous function $\pi(\theta)$ is a permissible prior for model $\mathcal{M} = \{p(\mathbf{x} \mid \theta), \mathbf{x} \in \mathcal{X}, \theta \in \Theta\}$ if:

1. for all $\mathbf{x} \in \mathcal{X}$, $\pi(\theta \mid \mathbf{x})$ is proper, that is, $\int_\Theta p(\mathbf{x} \mid \theta)\pi(\theta)\,d\theta < \infty$;
2. for some approximating compact sequence, the corresponding posterior sequence is expected logarithmically convergent to $\pi(\theta \mid \mathbf{x}) \propto p(\mathbf{x} \mid \theta)\pi(\theta)$.

The following theorem, whose proof is given in Appendix A, shows that, for one observation from a location model, the objective prior $\pi(\theta) = 1$ is permissible under mild conditions.

THEOREM 2. *Consider the model* $\mathcal{M} = \{f(x - \theta), \theta \in \mathbb{R}, x \in \mathbb{R}\}$, *where $f(t)$ is a density function on $\mathbb{R}$. If, for some $\varepsilon > 0$,*

$$\lim_{|t| \to 0} |t|^{1+\varepsilon} f(t) = 0, \tag{2.3}$$

*then $\pi(\theta) = 1$ is a permissible prior for the location model $\mathcal{M}$.*

EXAMPLE 3 (A nonpermissible constant prior in a location model). Consider the location model $\mathcal{M} \equiv \{p(x \mid \theta) = f(x - \theta), \theta \in \mathbb{R}, x > \theta + e\}$, where $f(t) = t^{-1}(\log t)^{-2}$, $t > e$. It is shown in Appendix B that, if $\pi(\theta) = 1$, then $\int_{\Theta_0} \kappa\{\pi(\theta \mid x) \mid \pi_0(\theta \mid x)\} p_0(x)\, dx = \infty$ for any compact set $\Theta_0 = [a, b]$ with $b - a \geq 1$; thus, $\pi(\theta) = 1$ is not a permissible prior for $\mathcal{M}$. Note that this model does not satisfy (2.3).

This is an interesting example because we are still dealing with a location density, so that $\pi(\theta) = 1$ is still the invariant (Haar) prior and, as such, satisfies numerous nice properties such as being exact frequentist matching (i.e., a Bayesian $100(1 - \alpha)\%$ credible set will also be a frequentist $100(1 - \alpha)\%$ confidence set; cf. equation (6.22) in Berger [2]). This is in stark contrast to the situation with the Fraser, Monette and Ng example. However, the basic fact remains that posteriors from uniform priors on large compact sets do not seem here to be well approximated (in terms of logarithmic divergence) by a uniform prior on the full parameter space. The suggestion is that this is a situation in which assessment of the "true" bounded parameter space is potentially needed.

Of course, a prior might be permissible for a larger sample size, even if it is not permissible for the minimal sample size. For instance, we suspect that $\pi(\theta) = 1$ is permissible for any location model having two or more independent observations.

The condition in the definition of permissibility that the posterior must be proper is not vacuous, as the following example shows.

EXAMPLE 4 (Mixture model). Let $\mathbf{x} = \{x_1, \ldots, x_n\}$ be a random sample from the mixture $p(x_i \mid \theta) = \frac{1}{2} \mathrm{N}(x \mid \theta, 1) + \frac{1}{2} \mathrm{N}(x \mid 0, 1)$, and consider the uniform prior function $\pi(\theta) = 1$. Since the likelihood function is bounded below by $2^{-n} \prod_{j=1}^{n} \mathrm{N}(x_j \mid 0, 1) > 0$, the integrated likelihood $\int_{-\infty}^{\infty} p(\mathbf{x} \mid \theta) \pi(\theta) \, d\theta = \int_{-\infty}^{\infty} p(\mathbf{x} \mid \theta) \, d\theta$ will diverge. Hence, the corresponding formal posterior is improper, and therefore the uniform prior is not a permissible prior function for this model. It can be shown that Jeffreys prior for this mixture model has the shape of an inverted bell, with a minimum value $1/2$ at $\mu = 0$; hence, it is also bounded from below and is, therefore, not a permissible prior for this model either.

Example 4 is noteworthy because it is very rare for the Jeffreys prior to yield an improper posterior in univariate problems. It is also of interest because there is no natural objective prior available for the problem. (There are data-dependent objective priors: see Wasserman [43].)

Theorem 2 can easily be modified to apply to models that can be transformed into a location model.

COROLLARY 1. *Consider* $\mathcal{M} \equiv \{p(x \mid \theta), \theta \in \Theta, x \in \mathcal{X}\}$. *If there are monotone functions* $y = y(x)$ *and* $\phi = \phi(\theta)$ *such that* $p(y \mid \phi) = f(y - \phi)$ *is a location model and there exists* $\varepsilon > 0$ *such that* $\lim_{|t| \to 0} |t|^{1+\varepsilon} f(t) = 0$, *then* $\pi(\theta) = |\phi'(\theta)|$ *is a permissible prior function for* $\mathcal{M}$.

The most frequent transformation is the log transformation, which converts a scale model into a location model. Indeed, this transformation yields the following direct analogue of Theorem 2.

COROLLARY 2. *Consider* $\mathcal{M} = \{p(x \mid \theta) = \theta^{-1} f(|x|/\theta), \theta > 0, x \in \mathbb{R}\}$, *a scale model where* $f(s)$, $s > 0$, *is a density function. If, for some* $\varepsilon > 0$,

$$\lim_{|t| \to \infty} |t|^{1+\varepsilon} e^t f(e^t) = 0, \tag{2.4}$$

*then* $\pi(\theta) = \theta^{-1}$ *is a permissible prior function for the scale model* $\mathcal{M}$.

EXAMPLE 5 (Exponential data). If $x$ is an observation from an exponential density, (2.4) becomes $|t|^{1+\varepsilon} e^t \exp(-e^t) \to 0$, as $|t| \to \infty$, which is true. From Corollary 2, $\pi(\theta) = \theta^{-1}$ is a permissible prior; indeed, $\pi_i(\theta) = (2i)^{-1} \theta^{-1}$, $e^{-i} \le \theta \le e^i$ is expected logarithmically convergent to $\pi(\theta)$.

EXAMPLE 6 (Uniform data). Let $x$ be one observation from the uniform distribution $\mathcal{M} = \{\mathrm{Un}(x \mid 0, \theta) = \theta^{-1}, x \in [0, \theta], \theta > 0\}$. This is a scale density, and equation (2.4) becomes $|t|^{1+\varepsilon} e^t \mathbf{1}_{\{0 < e^t < 1\}} \to 0$, as $|t| \to \infty$, which is indeed true. Thus, $\pi(\theta) = \theta^{-1}$ is a permissible prior function for $\mathcal{M}$.

The examples showing permissibility were for a single observation. Pleasantly, it is enough to establish permissibility for a single observation or, more generally, for the sample size necessary for posterior propriety of $\pi(\theta \mid \mathbf{x})$ because of the following theorem, which shows that expected logarithmic discrepancy is monotonically nonincreasing in sample size.

THEOREM 3 (Monotone expected logarithmic discrepancy). *Let $\mathcal{M} = \{p(\mathbf{x}_1, \mathbf{x}_2 \mid \theta) = p(\mathbf{x}_1 \mid \theta)p(\mathbf{x}_2 \mid \mathbf{x}_1, \theta), \mathbf{x}_1 \in \mathcal{X}_1, \mathbf{x}_2 \in \mathcal{X}_2, \theta \in \Theta\}$ be a parametric model. Consider a continuous improper prior $\pi(\theta)$ satisfying $m(\mathbf{x}_1) = \int_\Theta p(\mathbf{x}_1 \mid \theta)\pi(\theta)\, d\theta < \infty$ and $m(\mathbf{x}_1, \mathbf{x}_2) = \int_\Theta p(\mathbf{x}_1, \mathbf{x}_2 \mid \theta)\pi(\theta)\, d\theta < \infty$. For any compact set $\Theta_0 \subset \Theta$, let $\pi_0(\theta) = \pi(\theta)\mathbf{1}_{\Theta_0}(\theta)/\int_{\Theta_0} \pi(\theta)\, d\theta$. Then,*

$$
\begin{aligned}
\int\int_{\mathcal{X}_1 \times \mathcal{X}_2} &\kappa\{\pi(\cdot \mid \mathbf{x}_1, \mathbf{x}_2) \mid \pi_0(\cdot \mid \mathbf{x}_1, \mathbf{x}_2)\} m_0(\mathbf{x}_1, \mathbf{x}_2)\, d\mathbf{x}_1\, d\mathbf{x}_2 \\
&\leq \int_{\mathcal{X}_1} \kappa\{\pi(\cdot \mid \mathbf{x}_1) \mid \pi_0(\cdot \mid \mathbf{x}_1)\} m_0(\mathbf{x}_1)\, d\mathbf{x}_1,
\end{aligned}
$$
(2.5)

*where for $\theta \in \Theta_0$,*

$$\pi_0(\theta \mid \mathbf{x}_1, \mathbf{x}_2) = \frac{p(\mathbf{x}_1, \mathbf{x}_2 \mid \theta)\pi(\theta)}{m_0(\mathbf{x}_1, \mathbf{x}_2)},$$

$$m_0(\mathbf{x}_1, \mathbf{x}_2) = \int_{\Theta_0} p(\mathbf{x}_1, \mathbf{x}_2 \mid \theta)\pi(\theta)\, d\theta,$$

$$\pi_0(\theta \mid \mathbf{x}_1) = \frac{p(\mathbf{x}_1 \mid \theta)\pi(\theta)}{m_0(\mathbf{x}_1)},$$

$$m_0(\mathbf{x}_1) = \int_{\Theta_0} p(\mathbf{x}_1 \mid \theta)\pi(\theta)\, d\theta.$$

PROOF. The proof of this theorem is given in Appendix C. □

As an aside, the above result suggests that, as the sample size grows, the convergence of the posterior to normality given in Clarke [16] is monotone.

### 3. Reference priors.

3.1. *Definition of reference priors.* Key to the definition of reference priors is *Shannon expected information* (Shannon [38] and Lindley [36]).

DEFINITION 6 (Expected information). The information to be expected from one observation from model $\mathcal{M} \equiv \{p(\mathbf{x} \mid \theta), \mathbf{x} \in \mathcal{X}, \theta \in \Theta\}$, when the prior for $\theta$ is $q(\theta)$, is

$$I\{q \mid \mathcal{M}\} = \int\int_{\mathcal{X} \times \Theta} p(\mathbf{x} \mid \theta)q(\theta) \log \frac{p(\theta \mid \mathbf{x})}{q(\theta)}\, d\mathbf{x}\, d\theta$$

(3.1)
$$= \int_{\mathcal{X}} \kappa\{q \mid p(\cdot \mid \mathbf{x})\} p(\mathbf{x}) \, d\mathbf{x},$$

where $p(\theta \mid \mathbf{x}) = p(\mathbf{x} \mid \theta) q(\theta)/p(\mathbf{x})$ and $p(\mathbf{x}) = \int_\Theta p(\mathbf{x} \mid \theta) q(\theta) \, d\theta$.

Note that $\mathbf{x}$ here refers to the entire observation vector. It can have any dependency structure whatsoever (e.g., it could consist of $n$ normal random variables with mean zero, variance one and correlation $\theta$.) Thus, when we refer to a model henceforth, we mean the probability model for the actual complete observation vector. Although somewhat nonstandard, this convention is necessary here because reference prior theory requires the introduction of (artificial) independent replications of the entire experiment.

The amount of information $I\{q \mid \mathcal{M}\}$ to be expected from observing $\mathbf{x}$ from $\mathcal{M}$ depends on the prior $q(\theta)$: the sharper the prior the smaller the amount of information to be expected from the data. Consider now the information $I\{q \mid \mathcal{M}^k\}$ which may be expected from $k$ independent replications of $\mathcal{M}$. As $k \to \infty$, the sequence of realizations $\{\mathbf{x}_1, \ldots, \mathbf{x}_k\}$ would eventually provide any missing information about the value of $\theta$. Hence, as $k \to \infty$, $I\{q \mid \mathcal{M}^k\}$ provides a measure of the missing information about $\theta$ associated to the prior $q(\theta)$. Intuitively, a reference prior will be a permissible prior which maximizes the missing information about $\theta$ within the class $\mathcal{P}$ of priors compatible with any assumed knowledge about the value of $\theta$.

With a continuous parameter space, the missing information $I\{q \mid \mathcal{M}^k\}$ will typically diverge as $k \to \infty$, since an infinite amount of information would be required to learn the value of $\theta$. Likewise, the expected information is typically not defined on an unbounded set. These two difficulties are overcome with the following definition, that formalizes the heuristics described in Bernardo [10] and in Berger and Bernardo [7].

DEFINITION 7 [Maximizing Missing Information (MMI) Property]. Let $\mathcal{M} \equiv \{p(\mathbf{x} \mid \theta), \mathbf{x} \in \mathcal{X}, \theta \in \Theta \in \mathbb{R}\}$, be a model with one continuous parameter, and let $\mathcal{P}$ be a class of prior functions for $\theta$ for which $\int_\Theta p(\mathbf{x} \mid \theta) p(\theta) \, d\theta < \infty$. The function $\pi(\theta)$ is said to have the MMI property for model $\mathcal{M}$ given $\mathcal{P}$ if, for any compact set $\Theta_0 \in \Theta$ and any $p \in \mathcal{P}$,

(3.2)
$$\lim_{k \to \infty} \{I\{\pi_0 \mid \mathcal{M}^k\} - I\{p_0 \mid \mathcal{M}^k\}\} \geq 0,$$

where $\pi_0$ and $p_0$ are, respectively, the renormalized restrictions of $\pi(\theta)$ and $p(\theta)$ to $\Theta_0$.

The restriction of the definition to a compact set typically ensures the existence of the missing information for given $k$. That the missing information will diverge for large $k$ is handled by the device of simply insisting that the missing information for the reference prior be larger, as $k \to \infty$, than the missing information for any other candidate $p(\theta)$.

DEFINITION 8. A function $\pi(\theta) = \pi(\theta \mid \mathcal{M}, \mathcal{P})$ is a reference prior for model $\mathcal{M}$ given $\mathcal{P}$ if it is permissible and has the MMI property.

Implicit in this definition is that the reference prior on $\Theta$ will also be the reference prior on any compact subset $\Theta_0$. This is an attractive property that is often stated as the practical way to proceed when dealing with a restricted parameter space, but here it is simply a consequence of the definition.

Although we feel that a reference prior needs to be both permissible and have the MMI property, the MMI property is considerably more important. Thus, others have defined reference priors only in relation to this property, and Definition 7 is compatible with a number of these previous definitions in particular cases. Clarke and Barron [17] proved that, under appropriate regularity conditions, essentially those which guarantee asymptotic posterior normality, the prior which asymptotically maximizes the information to be expected by repeated sampling from $\mathcal{M} \equiv \{p(x \mid \theta), x \in \mathcal{X}, \theta \in \Theta \in \mathbb{R}\}$ is the Jeffreys prior,

$$(3.3) \quad \pi(\theta) = \sqrt{i(\theta)}, \qquad i(\theta) = -\int_{\mathcal{X}} p(x \mid \theta) \frac{\partial^2}{(\partial \theta)^2} \log[p(x \mid \theta)] \, dx$$

which, hence, is the reference prior under those conditions. Similarly, Ghosal and Samanta [27] gave conditions under which the prior, which asymptotically maximizes the information to be expected by repeated sampling from nonregular models of the form $\mathcal{M} \equiv \{p(x \mid \theta), x \in S(\theta), \theta \in \Theta \in \mathbb{R}\}$, where the support $S(\theta)$ is either monotonically decreasing or monotonically increasing in $\theta$, is

$$(3.4) \quad \pi(\theta) = \left| \int_{\mathcal{X}} p(x \mid \theta) \frac{\partial}{\partial \theta} \log[p(x \mid \theta)] \, dx \right|,$$

which is, therefore, the reference prior under those conditions.

3.2. *Properties of reference priors.* Some important properties of reference priors—generally regarded as required properties for any sensible procedure to derive objective priors—can be immediately deduced from their definition.

THEOREM 4 (Independence of sample size). *If data $\mathbf{x} = \{\mathbf{y}_1, \ldots, \mathbf{y}_n\}$ consists of a random sample of size $n$ from model $\mathcal{M} = \{p(\mathbf{y} \mid \theta), \mathbf{y} \in \mathcal{Y}, \theta \in \Theta\}$ with reference prior $\pi(\theta \mid \mathcal{M}, \mathcal{P})$, then $\pi(\theta \mid \mathcal{M}^n, \mathcal{P}) = \pi(\theta \mid \mathcal{M}, \mathcal{P})$, for any fixed sample size $n$.*

PROOF. This follows from the additivity of the information measure. Indeed, for any sample size $n$ and number of replicates $k$, $I\{q \mid \mathcal{M}^{nk}\} = nI\{q \mid \mathcal{M}^k\}$. $\square$

Note, however, that Theorem 4 requires **x** to be a random sample from the assumed model. If observations are dependent, as in time series or spatial models, the reference prior may well depend on the sample size (see, e.g., Berger and Yang [9] and Berger, de Oliveira and Sansó [8]).

THEOREM 5 (Compatibility with sufficient statistics). *Consider the model* $\mathcal{M} = \{p(\mathbf{x} \mid \theta), \mathbf{x} \in \mathcal{X}, \theta \in \Theta\}$ *with sufficient statistic* $\mathbf{t} = \mathbf{t}(\mathbf{x}) \in \mathcal{T}$, *and let* $\mathcal{M}_\mathbf{t} = \{p(\mathbf{t} \mid \theta), \mathbf{t} \in \mathcal{T}, \theta \in \Theta\}$ *be the corresponding model in terms of* **t**. *Then,* $\pi(\theta \mid \mathcal{M}, \mathcal{P}) = \pi(\theta \mid \mathcal{M}_\mathbf{t}, \mathcal{P})$.

PROOF. This follows because expected information is invariant under such transformation, so that, for all $k$, $I\{q \mid \mathcal{M}^k\} = I\{q \mid \mathcal{M}_\mathbf{t}^k\}$. □

THEOREM 6 (Consistency under reparametrization). *Consider the model* $\mathcal{M}_1 = \{p(\mathbf{x} \mid \theta), \mathbf{x} \in \mathcal{X}, \theta \in \Theta\}$, *let* $\phi(\theta)$ *be an invertible transformation of* $\theta$, *and let* $\mathcal{M}_1$ *be the model parametrized in terms of* $\phi$. *Then,* $\pi(\phi \mid \mathcal{M}_2, \mathcal{P})$ *is the prior density induced from* $\pi(\theta \mid \mathcal{M}_1, \mathcal{P})$ *by the appropriate probability transformation.*

PROOF. This follows immediately from the fact that the expected information is also invariant under one-to-one reparametrizations, so that, for all $k$, $I\{q_1 \mid \mathcal{M}_1^k\} = I\{q_2 \mid \mathcal{M}_2^k\}$, where $q_2(\phi) = q_1(\theta) \times |\partial\theta/\partial\phi|$. □

3.3. *Existence of the expected information.* The definition of a reference prior is clearly only useful if the $I\{\pi_0 \mid \mathcal{M}^k\}$ and $I\{p_0 \mid \mathcal{M}^k\}$ are finite for the (artificial) replications of $\mathcal{M}$. It is useful to write down conditions under which this will be so.

DEFINITION 9 (Standard prior functions). Let $\mathcal{P}_s$ be the class of strictly positive and continuous prior functions on $\Theta$ which have proper formal posterior distributions so that, when $p \in \mathcal{P}_s$,

$$(3.5) \qquad \forall \theta \in \Theta, \ p(\theta) > 0; \ \forall \mathbf{x} \in \mathcal{X} \qquad \int_\Theta p(\mathbf{x} \mid \theta) p(\theta) \, d\theta < \infty.$$

We call these the standard prior functions.

This will be the class of priors that we typically use to define the reference prior. The primary motivation for requiring a standard prior to be positive and continuous on $\Theta$ is that any prior not satisfying these conditions would not be accepted as being a reasonable candidate for an objective prior.

DEFINITION 10 (Standard models). Let $\mathcal{M} \equiv \{p(\mathbf{x} \mid \theta), \mathbf{x} \in \mathcal{X}, \theta \in \Theta \subset \mathbb{R}\}$ be a model with continuous parameter, and let $\mathbf{t}_k = \mathbf{t}_k(\mathbf{x}_1, \ldots, \mathbf{x}_k) \in \mathcal{T}_k$ be any sufficient statistic for the (artificial) $k$ replications of the experiment. ($\mathbf{t}_k$ could be just the observation vectors themselves.) The model $\mathcal{M}$ is said to be standard if, for any prior function $p(\theta) \in \mathcal{P}_s$ and any compact set $\Theta_0$,

$$\tag{3.6} I\{p_0 \mid \mathcal{M}^k\} < \infty,$$

where $p_0(\theta)$ is the proper prior obtained by restricting $p(\theta)$ to $\Theta_0$.

There are a variety of conditions under which satisfaction of (3.6) can be assured. Here is one of the simplest, useful when all $p(\mathbf{t}_k \mid \theta)$, for $\theta \in \Theta_0$, have the same support.

LEMMA 1. *For $p(\theta) \in \mathcal{P}_s$ and any compact set $\Theta_0$, (3.6) is satisfied if, for any $\theta \in \Theta_0$ and $\theta' \in \Theta_0$,*

$$\tag{3.7} \int_{\mathcal{T}_k} p(\mathbf{t}_k \mid \theta) \log \frac{p(\mathbf{t}_k \mid \theta)}{p(\mathbf{t}_k \mid \theta')} \, d\mathbf{t}_k < \infty.$$

PROOF. The proof is given in Appendix D. □

When the $p(\mathbf{t}_k \mid \theta)$ have different supports over $\theta \in \Theta_0$, the following lemma, whose proof is given in Appendix E, can be useful to verify (3.6).

LEMMA 2. *For $p(\theta) \in \mathcal{P}_s$ and any compact set $\Theta_0$, (3.6) is satisfied if:*

1. $H[p(\mathbf{t}_k \mid \theta)] \equiv -\int_{\mathcal{T}_k} p(\mathbf{t}_k \mid \theta) \log[p(\mathbf{t}_k \mid \theta)] \, d\mathbf{t}_k$ *is bounded below over $\Theta_0$.*
2. $\int_{\mathcal{T}_k} p_0(\mathbf{t}_k) \log[p_0(\mathbf{t}_k)] \, d\mathbf{t}_k > -\infty$, *where $p_0(\mathbf{t}_k)$ is the marginal likelihood from the uniform prior, that is, $p_0(\mathbf{t}_k) = L(\Theta_0)^{-1} \int_{\Theta_0} p(\mathbf{t}_k \mid \theta) \, d\theta$, with $L(\Theta_0)$ being the Lebesgue measure of $\Theta_0$.*

## 4. Determining the reference prior.

4.1. *An explicit expression for the reference prior.* Definition 8 does not provide a constructive procedure to derive reference priors. The following theorem provides an explicit expression for the reference prior, under certain mild conditions. Recall that $\mathbf{x}$ refers to the entire vector of observations from the model, while $\mathbf{x}^{(k)} = (\mathbf{x}_1, \ldots, \mathbf{x}_k)$ refers to a vector of (artificial) independent replicates of these vector observations from the model. Finally, let $\mathbf{t}_k = \mathbf{t}_k(\mathbf{x}_1, \ldots, \mathbf{x}_k) \in \mathcal{T}_k$ be any sufficient statistic for the replicated observations. While $\mathbf{t}_k$ could just be $\mathbf{x}^{(k)}$ itself, it is computationally convenient to work with sufficient statistics if they are available.

THEOREM 7 (*Explicit form of the reference prior*). *Assume a standard model* $\mathcal{M} \equiv \{p(\mathbf{x} \mid \theta), \mathbf{x} \in \mathcal{X}, \theta \in \Theta \subset \mathbb{R}\}$ *and the standard class* $\mathcal{P}_s$ *of candidate priors. Let* $\pi^*(\theta)$ *be a continuous strictly positive function such that the corresponding formal posterior*

$$(4.1) \qquad \pi^*(\theta \mid \mathbf{t}_k) = \frac{p(\mathbf{t}_k \mid \theta)\pi^*(\theta)}{\int_\Theta p(\mathbf{t}_k \mid \theta)\pi^*(\theta)\,d\theta}$$

*is proper and asymptotically consistent (see Appendix F), and define, for any interior point $\theta_0$ of $\Theta$,*

$$(4.2) \qquad f_k(\theta) = \exp\left\{\int_{\mathcal{T}_k} p(\mathbf{t}_k \mid \theta)\log[\pi^*(\theta \mid \mathbf{t}_k)]\,d\mathbf{t}_k\right\} \quad \text{and}$$

$$(4.3) \qquad f(\theta) = \lim_{k \to \infty} \frac{f_k(\theta)}{f_k(\theta_0)}.$$

*If* (i) *each $f_k(\theta)$ is continuous and, for any fixed $\theta$ and sufficiently large $k$, $\{f_k^0(\theta)/f_k^0(\theta_0)\}$ is either monotonic in $k$ or is bounded above by some $h(\theta)$ which is integrable on any compact set, and* (ii) *$f(\theta)$ is a permissible prior function, then $\pi(\theta \mid \mathcal{M}, \mathcal{P}_s) = f(\theta)$ is a reference prior for model $\mathcal{M}$ and prior class $\mathcal{P}_s$.*

PROOF. The proof of this theorem is given in Appendix F. □

Note that the choice of $\pi^*$ is essentially arbitrary and, hence, can be chosen for computational convenience. Also, the choice of $\theta_0$ is immaterial. Finally, note that no compact set is mentioned in the theorem; that is, the defined reference prior works simultaneously for all compact subsets of $\Theta$.

EXAMPLE 7 (Location model). To allow for the dependent case, we write the location model for data $\mathbf{x} = (x_1, \ldots, x_n)$ as $f(x_1 - \theta, \ldots, x_n - \theta)$, where we assume $\Theta = \mathbb{R}$. To apply Theorem 7, choose $\pi^*(\theta) = 1$. Then, because of the translation invariance of the problem, it is straightforward to show that (4.2) reduces to $f_k(\theta) = c_k$, not depending on $\theta$. It is immediate from (4.3) that $f(\theta) = 1$, and condition (a) of Theorem 7 is also trivially satisfied. [Note that this is an example of choosing $\pi^*(\theta)$ conveniently; any other choice would have resulted in a much more difficult analysis.]

It follows that, if the model is a standard model and $\pi(\theta) = 1$ is permissible for the model [certainly satisfied if (2.3) holds], then $\pi(\theta) = 1$ is the reference prior among the class of all standard priors. Note that there is additional work that is needed to verify that the model is a standard model. Easiest is typically to verify (3.7), which is easy for most location families.

It is interesting that no knowledge of the asymptotic distribution of the posterior is needed for this result. Thus, the conclusion applies equally to the normal distribution and to the distribution with density $f(x - \theta) = \exp(x - \theta)$, for $x > \theta$, which is not asymptotically normal.

The key feature making the analysis for the location model so simple was that (4.2) was a constant. A similar result will hold for any suitably invariant statistical model if $\pi^*(\theta)$ is chosen to be the Haar density (or right-Haar density in multivariable models); then, (4.2) becomes a constant times the right-Haar prior. For instance, scale-parameter problems can be handled in this way, although one can, of course, simply transform them to a location parameter problem and apply the above result. For a scale-parameter problem, the reference prior is, of course, $\pi(\theta) = \theta^{-1}$.

EXAMPLE 8. A model for which nothing is known about reference priors is the uniform model with support on $(a_1(\theta), a_2(\theta))$,

$$(4.4) \quad \mathcal{M} = \left\{ \mathrm{Un}(x \mid a_1(\theta), a_2(\theta)) = \frac{1}{a_2(\theta) - a_1(\theta)}, a_1(\theta) < x < a_2(\theta) \right\},$$

where $\theta > \theta_0$ and $0 < a_1(\theta) < a_2(\theta)$ are both strictly monotonic increasing functions on $\Theta = (\theta_0, \infty)$ with derivatives satisfying $0 < a_1'(\theta) < a_2'(\theta)$. This is not a regular model, has no group invariance structure and does not belong to the class of nonregular models analyzed in Ghosal and Samanta [27]. The following theorem gives the reference prior for the model (4.4). Its proof is given in Appendix G.

THEOREM 8. *Consider the model (4.4). Define*

$$(4.5) \qquad b_j \equiv b_j(\theta) = \frac{a_2'(\theta) - a_1'(\theta)}{a_j'(\theta)}, \qquad j = 1, 2.$$

*Then the reference prior of $\theta$ for the model (4.4) is*

$$\pi(\theta) = \frac{a_2'(\theta) - a_1'(\theta)}{a_2(\theta) - a_1(\theta)}$$

$$(4.6) \qquad \times \exp\left\{ b_1 + \frac{1}{b_1 - b_2}\left[ b_1 \psi\left(\frac{1}{b_1}\right) - b_2 \psi\left(\frac{1}{b_2}\right)\right]\right\},$$

*where $\psi(z)$ is the digamma function defined by $\psi(z) = \frac{d}{dz}\log(\Gamma(z))$ for $z > 0$.*

EXAMPLE 9 [Uniform distribution on $(\theta, \theta^2), \theta > 1$]. This is a special case of Theorem 8, with $\theta_0 = 1$, $a_1(\theta) = \theta$ and $a_2(\theta) = \theta^2$. Then, $b_1 = 2\theta - 1$ and $b_2 = (2\theta - 1)/(2\theta)$. It is easy to show that $b_2^{-1} = b_1^{-1} + 1$. For the digamma function (see Boros and Moll [14]), $\psi(z+1) = \psi(z) + 1/z$, for $z > 0$, so that $\psi(1/b_1) = \psi(1/b_2) - b_1$. The reference prior (4.6) thus becomes

$$\pi(\theta) = \frac{2\theta - 1}{\theta(\theta - 1)} \exp\left\{ b_1 + \frac{1}{b_1 - b_2}\left[ b_1 \psi\left(\frac{1}{b_2}\right) - b_1^2 - b_2 \psi\left(\frac{1}{b_2}\right)\right]\right\}$$

$$(4.7) \qquad = \frac{2\theta - 1}{\theta(\theta - 1)} \exp\left\{ -\frac{b_1 b_2}{b_1 - b_2} + \psi\left(\frac{1}{b_2}\right)\right\}$$

$$\propto \frac{2\theta-1}{\theta(\theta-1)}\exp\left\{\psi\left(\frac{2\theta}{2\theta-1}\right)\right\},$$

the last equation following from the identity $b_1 b_2/(b_1 - b_2) = b_2^{-1} - b_1^{-1} = 1$.

4.2. *Numerical computation of the reference prior.* Analytical derivation of reference priors may be technically demanding in complex models. However, Theorem 7 may also be used to obtain an approximation to the reference prior through numerical evaluation of equation (4.2). Moderate values of $k$ (to simulate the asymptotic posterior) will often yield a good approximation to the reference prior. The appropriate pseudo code is:

ALGORITHM.

1. Starting values:
   choose a moderate value for $k$;
   choose an arbitrary positive function $\pi^*(\theta)$, say $\pi^*(\theta) = 1$;
   choose the number $m$ of samples to be simulated.
2. For any given $\theta$ value, **repeat**, for $j = 1, \ldots, m$:
   simulate a random sample $\{\mathbf{x}_{1j}, \ldots, \mathbf{x}_{kj}\}$ of size $k$ from $p(\mathbf{x} \mid \theta)$;
   compute numerically the integral $c_j = \int_\Theta \prod_{i=1}^k p(\mathbf{x}_{ij} \mid \theta) \pi^*(\theta) \, d\theta$;
   evaluate $r_j(\theta) = \log[\prod_{i=1}^k p(\mathbf{x}_{ij} \mid \theta) \pi^*(\theta)/c_j]$.
3. Compute $\pi(\theta) = \exp[m^{-1} \sum_{j=1}^m r_j(\theta)]$ and **store** the pair $\{\theta, \pi(\theta)\}$.
4. **Repeat** routines (2) and (3) for all $\theta$ values for which the pair $\{\theta, \pi(\theta)\}$ is required.

If desired, a continuous approximation to $\pi(\theta)$ may easily be obtained from the computed points using standard interpolation techniques.

We first illustrate the computation in an example for which the reference prior is known to enable comparison of the numerical accuracy of the approximation.

EXAMPLE 10 [Uniform distribution on $(\theta, \theta^2)$, continued]. Consider again the uniform distribution on $(\theta, \theta^2)$ discussed in Example 9, where the reference prior was analytically given in (4.7).

Figure 1 presents the reference prior numerically calculated with the algorithm for nine $\theta$ values, uniformly log-spaced and rescaled to have $\pi(2) = 1$; $m = 1000$ samples of $k = 500$ observations were used to compute each of the nine $\{\theta_i, \pi(\theta_i)\}$ points. These nine points are clearly almost perfectly fitted by the exact reference prior (4.7), shown by a continuous line; indeed, the nine points were accurate to within four decimal points.

This numerical computation was done before the analytic reference prior was obtained for the problem, and a nearly perfect fit to the nine $\theta$ values was obtained by the function $\pi(\theta) = 1/(\theta - 1)$, which was thus guessed to be the actual reference prior. This guess was wrong, but note that (4.7) over the computed range is indeed nearly proportional to $1/(\theta - 1)$.

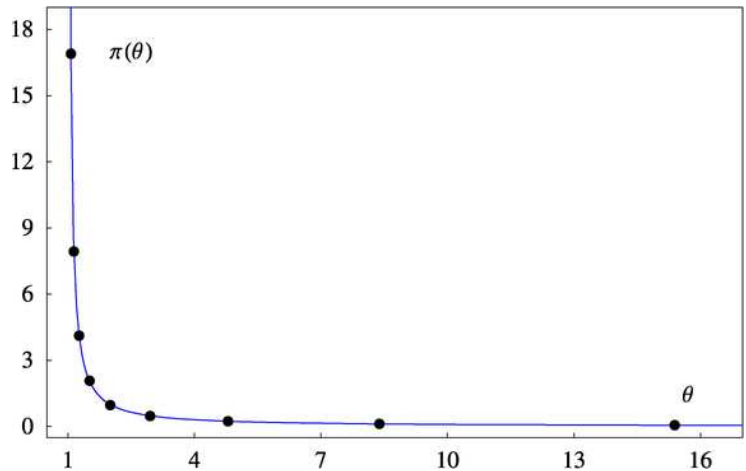

FIG. 1. *Numerical reference prior for the uniform model on* $[\theta, \theta^2]$.

We now consider an example for which the reference prior is not known and, indeed, appears to be extremely difficult to determine analytically.

EXAMPLE 11 (Triangular distribution). The use of a symmetric triangular distribution on $(0, 1)$ can be traced back to the 18th century to Simpson [39]. Schmidt [37] noticed that this pdf is the density of the mean of two i.i.d. uniform random variables on the interval $(0, 1)$.

The nonsymmetric standard triangular distribution on $(0, 1)$,

$$p(x \mid \theta) = \begin{cases} 2x/\theta, & \text{for } 0 < x \leq \theta, \\ 2(1-x)/(1-\theta), & \text{for } \theta < x < 1, \end{cases} \qquad 0 < \theta < 1,$$

was first studied by Ayyangar [1]. Johnson and Kotz [33] revisited nonsymmetric triangular distributions in the context of modeling prices. The triangular density has a unique mode at $\theta$ and satisfies $\Pr[x \leq \theta] = \theta$, a property that can be used to obtain an estimate of $\theta$ based on the empirical distribution function. The nonsymmetric triangular distribution does not possess a useful reduced sufficient statistic. Also, although $\log[p(x \mid \theta)]$ is differentiable for all $\theta$ values, the formal Fisher information function is strictly negative, so Jeffreys prior does not exist.

Figure 2 presents a numerical calculation of the reference prior at thirteen $\theta$ values, uniformly spaced on $(0, 1)$ and rescaled to have $\pi(1/2) = 2/\pi$; $m = 2500$ samples of $k = 2000$ observations were used to compute each of the thirteen $\{\theta_i, \pi(\theta_i)\}$ points. Interestingly, these points are nearly perfectly fitted by the (proper) prior $\pi(\theta) = \text{Be}(\theta \mid 1/2, 1/2) \propto \theta^{-1/2}(1-\theta)^{-1/2}$, shown by a continuous line.

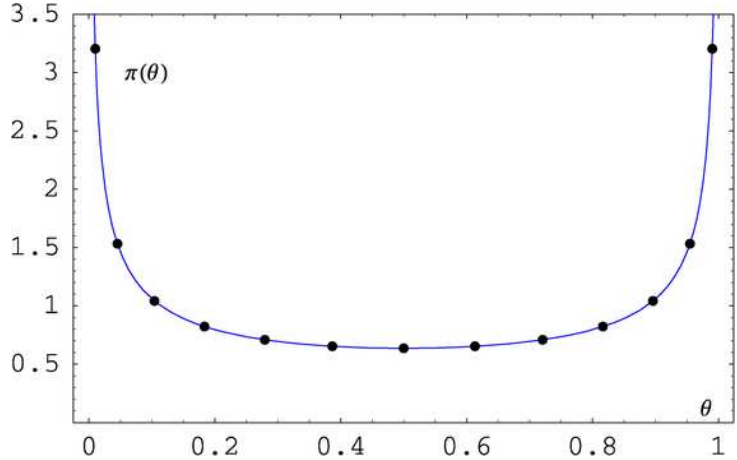

Fig. 2. *Numerical reference prior for the triangular model.*

Analytical derivation of the reference prior does not seem to be feasible in this example, but there is an interesting heuristic argument which suggests that the $\text{Be}(\theta \mid 1/2, 1/2)$ prior is indeed the reference prior for the problem. The argument begins by noting that, if $\tilde{\theta}_k$ is a consistent, asymptotically sufficient estimator of $\theta$, one would expect that, for large $k$,

$$\int_{\mathcal{T}} p(\mathbf{t}_k \mid \theta) \log[\pi_0(\theta \mid \mathbf{t}_k)] \, d\mathbf{t}_k \approx \int_{\mathcal{T}} p(\tilde{\theta}_k \mid \theta) \log[\pi_0(\theta \mid \tilde{\theta}_k)] \, d\tilde{\theta}_k$$
$$\approx \log[\pi_0(\theta \mid \tilde{\theta}_k)]|_{\tilde{\theta}_k = \theta},$$

since the sampling distribution of $\tilde{\theta}_k$ will concentrate on $\theta$. Thus, using (4.2) and (4.3), the reference prior should be

$$(4.8) \qquad \pi(\theta) = \pi_0(\theta \mid \tilde{\theta}_k)|_{\tilde{\theta}_k = \theta} \propto p(\tilde{\theta}_k \mid \theta)|_{\tilde{\theta}_k = \theta}.$$

For the triangular distribution, a consistent estimator of $\theta$ can be obtained as the solution to the equation $F_k(t) = t$, where $F_k(t)$ is the empirical distribution function corresponding to a random sample of size $k$. Furthermore, one can show that this solution, $\theta_k^*$, is asymptotically normal $\text{N}(\theta_k^* \mid \theta, s(\theta)/\sqrt{k})$, where $s(\theta) = \sqrt{\theta(1-\theta)}$. Plugging this into (4.8) would yield the $\text{Be}(\theta \mid 1/2, 1/2)$ prior as the reference prior. To make this argument precise, of course, one would have to verify that the above heuristic argument holds and that $\theta_k^*$ is asymptotically sufficient.

**5. Conclusions and generalizations.** The formalization of the notions of permissibility and the MMI property—the two keys to defining a reference prior—are of interest in their own right, but happened to be a by-product of the main goal, which was to obtain the explicit representation of a reference prior given in Theorem 7. Because of this explicit representation and, as illustrated in the examples following the theorem, one can:

- Have a single expression for calculating the reference prior, regardless of the asymptotic nature of the posterior distribution.
- Avoid the need to do computations over approximating compact parameter spaces.
- Develop a fairly simple numerical technique for computing the reference prior in situations where analytic determination is too difficult.
- Have, as immediate, the result that the reference prior on any compact subset of the parameter space is simply the overall reference prior constrained to that set.

The main limitation of the paper is the restriction to single parameter models. It would obviously be very useful to be able to generalize the results to deal with nuisance parameters.

The results concerning permissibility essentially generalize immediately to the multi-parameter case. The MMI property (and hence formal definition of a reference prior) can also be generalized to the multi-parameter case, following Berger and Bernardo [7] (although note that there were heuristic elements to that generalization). The main motivation for this paper, however, was the explicit representation for the reference prior that was given in Theorem 7, and, unfortunately, there does not appear to be an analogue of this explicit representation in the multi-parameter case. Indeed, we have found that any generalizations seem to require expressions that involve limits over approximating compact sets, precisely the feature of reference prior computation that we were seeking to avoid.

## APPENDIX A: PROOF OF THEOREM 2

By the invariance of the model, $p(x) = \int_\Theta f(x - \theta) \pi(\theta) \, d\theta = 1$ and $\pi(\theta \mid x) = f(x - \theta)$. To verify (ii) of Definition 5, choose $\Theta_i = [-i, i]$. Then $\pi_i(\theta \mid x) = f(x - \theta)/[2i p_i(x)]$, $\theta \in \Theta_i$, where

$$p_i(x) = \frac{1}{2i} \int_{-i}^{i} f(x - \theta) \, d\theta = \frac{1}{2i}(F(x + i) - F(x - i)),$$

with $F(x) = \int_{-\infty}^{x} f(t) \, dt$. The logarithmic discrepancy between $\pi_i(\theta \mid x)$ and $\pi(\theta \mid x)$ is

$$\kappa\{\pi(\cdot \mid x) \mid \pi_i(\cdot \mid x)\} = \int_{-i}^{i} \pi_i(\theta \mid x) \log \frac{\pi_i(\theta \mid x)}{\pi(\theta \mid x)} \, d\theta$$

$$= \int_{-i}^{i} \pi_i(\theta \mid x) \log \frac{1}{2ip_i(x)} d\theta$$

$$= -\log[F(x+i) - F(x-i)],$$

and the expected discrepancy is

$$\int_{-\infty}^{\infty} \kappa\{\pi(\cdot \mid x) \mid \pi_i(\cdot \mid x)\} p_i(x) \, dx$$

$$= -\frac{1}{2i} \int_{-\infty}^{\infty} [F(x+i) - F(x-i)] \log[F(x+i) - F(x-i)] \, dx$$

$$= \int_{-\infty}^{-4} g(y,i) \, dy + \int_{-4}^{2} g(y,i) \, dy + \int_{2}^{\infty} g(y,i) \, dy = J_1 + J_2 + J_3,$$

where, using the transformation $y = (x-i)/i$,

$$g(y,i) = -\{F[(y+2)i] - F(yi)\} \log\{F[(y+2)i] - F(yi)\}.$$

Notice that for fixed $y \in (-4, 2)$, as $i \to \infty$,

$$F[(y+2)i] - F(yi) \to \begin{cases} 0, & \text{if } y \in (-4, -2), \\ 1, & \text{if } y \in (-2, 0), \\ 0, & \text{if } y \in (0, 2). \end{cases}$$

Since $-v \log v \leq e^{-1}$ for $0 \leq v \leq 1$, the dominated convergence theorem can be applied to $J_2$, so that $J_2$ converges to 0 as $i \to \infty$. Next, when $i$ is large enough and, for any $y \geq 2$,

$$F[(y+2)i] - F(yi) \leq \int_{yi}^{(y+2)i} \frac{1}{t^{1+\varepsilon}} dt = \frac{1}{\varepsilon}\left(\frac{1}{(yi)^\varepsilon} - \frac{1}{[(y+2)i]^\varepsilon}\right)$$

$$= \frac{(1+2/y)^\varepsilon - 1}{\varepsilon i^\varepsilon (y+2)^\varepsilon} \leq \frac{2^\varepsilon}{i^\varepsilon y(y+2)^\varepsilon}$$

the last inequality holding since, for $0 \leq v \leq 1$, $(1+v)^\varepsilon - 1 \leq \varepsilon 2^{\varepsilon-1} v$. Using the fact that $-v \log v$ is monotone increasing in $0 \leq v \leq e^{-1}$, we have

$$J_3 \leq -\frac{1}{2} \int_{2}^{\infty} \frac{2^\varepsilon}{i^\varepsilon y(y+2)^\varepsilon} \log \frac{2^\varepsilon}{i^\varepsilon y(y+2)^\varepsilon} \, dy,$$

which converges to 0 as $i \to \infty$. It may similarly be shown that $J_1$ converges to 0 as $i \to \infty$. Consequently, $\{\pi_i(\theta \mid x)\}_{i=1}^{\infty}$ is expected logarithmically convergent to $\pi(\theta \mid x)$, and thus, $\pi(\theta) = 1$ is permissible.

## APPENDIX B: DERIVATION OF RESULTS IN EXAMPLE 3

Consider a location family, $p(x \mid \theta) = f(x - \theta)$, where $x \in \mathbb{R}$ and $\theta \in \Theta = \mathbb{R}$, and $f$ is given by $f(x) = x^{-1}(\log x)^{-2} 1_{(e,\infty)}(x)$. Choose $\pi(\theta) = 1$ and $\Theta_0 =$

$[a, b]$ such that $L \equiv b - a \geq 1$. Then,

$$Lp_0(x) = \int_a^b f(x - \theta)\, d\theta$$

$$= \begin{cases} \dfrac{1}{\log(-b-x)} - \dfrac{1}{\log(-a-x)}, & \text{if } x \leq -b - e, \\ 1 - \dfrac{1}{\log(-a-x)}, & \text{if } -b - e < x \leq -a - e, \\ 0, & \text{if } x > -a - e. \end{cases}$$

The logarithmic discrepancy between $\pi_0(\theta \mid x)$ and $\pi(\theta \mid x)$ is

$$\kappa\{\pi(\cdot \mid x) \mid \pi_0(\cdot \mid x)\} = \int_a^b \pi_0(\theta \mid x) \log \frac{\pi_0(\theta \mid x)}{\pi(\theta \mid x)}\, d\theta$$

$$= \int_a^b \pi_0(\theta \mid x) \log \frac{1}{Lp_0(x)}\, d\theta = -\log[Lp_0(x)].$$

Then the expected discrepancy is

$$E_0^x \kappa\{\pi(\cdot \mid x) \mid \pi_0(\cdot \mid x)\}$$

$$\equiv \int_{-\infty}^{\infty} p_0(x) \kappa\{\pi(\cdot \mid x) \mid \pi_0(\cdot \mid x)\}\, dx$$

$$= -\frac{1}{L} \int_{-\infty}^{\infty} Lp_0(x) \log[Lp_0(x)]\, dx$$

$$\geq -\frac{1}{L} \int_{-\infty}^{-b-e} \left\{ \frac{1}{\log(-b-x)} - \frac{1}{\log(-a-x)} \right\}$$

$$\times \log\left\{ \frac{1}{\log(-b-x)} - \frac{1}{\log(-a-x)} \right\} dx$$

$$= -\frac{1}{L} \int_e^{\infty} \left\{ \frac{1}{\log(t)} - \frac{1}{\log(t+L)} \right\} \log\left\{ \frac{1}{\log(t)} - \frac{1}{\log(t+L)} \right\} dt$$

$$\geq -\frac{1}{L} \int_{Le}^{\infty} \left\{ \int_t^{t+L} \frac{1}{x \log^2(x)}\, dx \right\} \log\left\{ \int_t^{t+L} \frac{1}{x \log^2(x)}\, dx \right\} dt.$$

Making the transformation $y = t/L$, the right-hand side equals

$$-\int_e^{\infty} g_L(y) \log\{g_L(y)\}\, dy,$$

where

$$g_L(y) = \int_{yL}^{(y+1)L} \frac{1}{x(\log x)^2}\, dx = \frac{1}{\log(yL)} - \frac{1}{\log((y+1)L)}$$

$$= \frac{1}{\log(y) + \log(L)} - \frac{1}{\log(y) + \log(1 + 1/y) + \log(L)}$$

$$= \frac{\log(1 + 1/y)}{[\log(y) + \log(L)][\log(y+1) + \log(L)]}.$$

Because $\log(1 + 1/y) > 1/(y+1)$, for $y \geq e$,

$$g_L(y) \geq \frac{1}{(y+1)[\log(y+1) + \log(L)]^2}.$$

Since $-p\log(p)$ is an increasing function of $p \in (0, e^{-1})$, it follows that

$$E_0^x \kappa\{\pi(\cdot \mid x) \mid \pi_0(\cdot \mid x)\} \geq J_1 + J_2,$$

where

$$J_1 = \int_e^\infty \frac{\log(y+1)}{(y+1)[\log(y+1) + \log(L)]^2} \, dy,$$

$$J_2 = \int_e^\infty \frac{2\log[\log(y+1) + \log(L)]}{(y+1)[\log(y+1) + \log(L)]^2} \, dy.$$

Clearly $J_1 = \infty$ and $J_2$ is finite, so $E_0^x \kappa\{\pi(\cdot \mid x) \mid \pi_0(\cdot \mid x)\}$ does not exist.

## APPENDIX C: PROOF OF THEOREM 3

First,

$$\int_{\mathcal{X}_1 \times \mathcal{X}_2} \kappa\{\pi(\cdot \mid \mathbf{x}_1, \mathbf{x}_2) \mid \pi_0(\cdot \mid \mathbf{x}_1, \mathbf{x}_2)\} m_0(\mathbf{x}_1, \mathbf{x}_2) \, d\mathbf{x}_1 \, d\mathbf{x}_2$$

$$= \int_{\mathcal{X}_1 \times \mathcal{X}_2} \int_{\Theta_0} \log\left\{\frac{\pi_0(\theta \mid \mathbf{x}_1, \mathbf{x}_2)}{\pi(\theta \mid \mathbf{x}_1, \mathbf{x}_2)}\right\} \pi_0(\theta) p(\mathbf{x}_1, \mathbf{x}_2 \mid \theta) \, d\theta \, d\mathbf{x}_1 \, d\mathbf{x}_2$$

$$= \int_{\mathcal{X}_1 \times \mathcal{X}_2} \int_{\Theta_0} \log\left\{\frac{\pi_0(\theta) m(\mathbf{x}_1, \mathbf{x}_2)}{\pi(\theta) m_0(\mathbf{x}_1, \mathbf{x}_2)}\right\} \pi_0(\theta) p(\mathbf{x}_1, \mathbf{x}_2 \mid \theta) \, d\theta \, d\mathbf{x}_1 \, d\mathbf{x}_2$$

$$= \int_{\Theta_0} \log\left\{\frac{\pi_0(\theta)}{\pi(\theta)}\right\} \pi_0(\theta) \, d\theta$$

(C.1)
$$+ \int_{\mathcal{X}_1 \times \mathcal{X}_2} \log\left\{\frac{m(\mathbf{x}_1, \mathbf{x}_2)}{m_0(\mathbf{x}_1, \mathbf{x}_2)}\right\} m_0(\mathbf{x}_1, \mathbf{x}_2) \, d\mathbf{x}_1 \, d\mathbf{x}_2$$

$$= J_0 + \int_{\mathcal{X}_1} \int_{\mathcal{X}_2} \log\left\{\frac{m(\mathbf{x}_2 \mid \mathbf{x}_1) m(\mathbf{x}_1)}{m_0(\mathbf{x}_2 \mid \mathbf{x}_1) m_0(\mathbf{x}_1)}\right\}$$

$$\times m_0(\mathbf{x}_2 \mid \mathbf{x}_1) m_0(\mathbf{x}_1) \, d\mathbf{x}_1 \, d\mathbf{x}_2$$

$$\equiv J_0 + J_1 + J_2,$$

where $J_0 = \int_{\Theta_0} \log\{\pi_0(\theta)/\pi(\theta)\} \pi_0(\theta) \, d\theta$,

$$J_1 = \int_{\mathcal{X}_1} \log\left\{\frac{m(\mathbf{x}_1)}{m_0(\mathbf{x}_1)}\right\} m_0(\mathbf{x}_1) \, d\mathbf{x}_1,$$

$$J_2 = \int_{\mathcal{X}_1} \left(\int_{\mathcal{X}_2} \log\left\{\frac{m(\mathbf{x}_2 \mid \mathbf{x}_1)}{m_0(\mathbf{x}_2 \mid \mathbf{x}_1)}\right\} m_0(\mathbf{x}_2 \mid \mathbf{x}_1) \, d\mathbf{x}_2\right) m_0(\mathbf{x}_1) \, d\mathbf{x}_1.$$

By assumption, $J_0$ is finite. Note that both $m_0(\mathbf{x}_2 \mid \mathbf{x}_1)$ and $m(\mathbf{x}_2 \mid \mathbf{x}_1) = m(\mathbf{x}_1, \mathbf{x}_2)/m(\mathbf{x}_1) = \int_\Theta p(\mathbf{x}_2 \mid \mathbf{x}_1, \theta)\pi(\theta)\,d\theta$ are proper densities. Because $\log(t)$ is concave on $(0, \infty)$, we have

$$J_2 \leq \int_{\mathcal{X}_1} \log\left\{\int_{\mathcal{X}_2} \frac{m(\mathbf{x}_2 \mid \mathbf{x}_1)}{m_0(\mathbf{x}_2 \mid \mathbf{x}_1)} m_0(\mathbf{x}_2 \mid \mathbf{x}_1)\,d\mathbf{x}_2\right\} m_0(\mathbf{x}_1)\,d\mathbf{x}_1 = 0.$$

By the same argument leading to (C.1), one can show that

$$\int_{\mathcal{X}_1} \kappa\{\pi(\cdot \mid \mathbf{x}_1) \mid \pi_0(\cdot \mid \mathbf{x}_1)\} m_0(\mathbf{x}_1)\,d\mathbf{x}_1 = J_0 + J_1.$$

The result is immediate.

## APPENDIX D: PROOF OF LEMMA 1

Clearly

$$I\{p_0 \mid \mathcal{M}^k\} \equiv \int_{\Theta_0} p_0(\theta) \int_{\mathcal{T}_k} p(\mathbf{t}_k \mid \theta) \log\left[\frac{p_0(\theta \mid \mathbf{t}_k)}{p_0(\theta)}\right] d\mathbf{t}_k\,d\theta$$

$$= \int_{\Theta_0} p_0(\theta) \int_{\mathcal{T}_k} p(\mathbf{t}_k \mid \theta) \log\left[\frac{p(\mathbf{t}_k \mid \theta)}{p_0(\mathbf{t}_k)}\right] d\mathbf{t}_k\,d\theta$$

$$\leq \sup_{\theta \in \Theta_0} \int_{\mathcal{T}_k} p(\mathbf{t}_k \mid \theta) \log\left[\frac{p(\mathbf{t}_k \mid \theta)}{p_0(\mathbf{t}_k)}\right] d\mathbf{t}_k.$$

Writing $p_0(\mathbf{t}_k) = \int_{\Theta_0} p(\mathbf{t}_k \mid \theta')p_0(\theta')\,d\theta'$, note by convexity of $[-\log]$ that

$$\int_{\mathcal{T}_k} p(\mathbf{t}_k \mid \theta) \log\left[\frac{p(\mathbf{t}_k \mid \theta)}{p_0(\mathbf{t}_k)}\right] d\mathbf{t}_k$$

$$= -\int_{\mathcal{T}_k} p(\mathbf{t}_k \mid \theta) \log\left[\int_{\Theta_0} \frac{p(\mathbf{t}_k \mid \theta')}{p(\mathbf{t}_k \mid \theta)} p_0(\theta')\,d\theta'\right] d\mathbf{t}_k$$

$$\leq -\int_{\mathcal{T}_k} p(\mathbf{t}_k \mid \theta) \left[\int_{\Theta_0} \log\left[\frac{p(\mathbf{t}_k \mid \theta')}{p(\mathbf{t}_k \mid \theta)}\right] p_0(\theta')\,d\theta'\right] d\mathbf{t}_k$$

$$= -\int_{\Theta_0} \int_{\mathcal{T}_k} p(\mathbf{t}_k \mid \theta) \log\left[\frac{p(\mathbf{t}_k \mid \theta')}{p(\mathbf{t}_k \mid \theta)}\right] d\mathbf{t}_k\, p_0(\theta')\,d\theta'$$

$$\leq -\inf_{\theta' \in \Theta_0} \int_{\mathcal{T}_k} p(\mathbf{t}_k \mid \theta) \log\left[\frac{p(\mathbf{t}_k \mid \theta')}{p(\mathbf{t}_k \mid \theta)}\right] d\mathbf{t}_k.$$

Combining this with (D.1) yields

$$I\{p_0 \mid \mathcal{M}^k\} \leq \sup_{\theta \in \Theta_0} \sup_{\theta' \in \Theta_0} \int_{\mathcal{T}_k} p(\mathbf{t}_k \mid \theta) \log\left[\frac{p(\mathbf{t}_k \mid \theta)}{p(\mathbf{t}_k \mid \theta')}\right] d\mathbf{t}_k,$$

from which the result follows.

## APPENDIX E: PROOF OF LEMMA 2

Let $p_0(\theta \mid \mathbf{t}_k)$ be the posterior of $\theta$ under $p_0$, that is, $p(\mathbf{t}_k \mid \theta) p_0(\theta) / p_0(\mathbf{t}_k)$. Note that

$$I\{p_0 \mid \mathcal{M}^k\} \equiv \int_{\Theta_0} p_0(\theta) \int_{\mathcal{T}_k} p(\mathbf{t}_k \mid \theta) \log\left[\frac{p_0(\theta \mid \mathbf{t}_k)}{p_0(\theta)}\right] d\mathbf{t}_k \, d\theta$$

$$= \int_{\Theta_0} \int_{\mathcal{T}_k} p_0(\theta) p(\mathbf{t}_k \mid \theta) \log\left[\frac{p(\mathbf{t}_k \mid \theta)}{p_0(\mathbf{t}_k)}\right] d\mathbf{t}_k \, d\theta$$

$$= \int_{\Theta_0} p_0(\theta) \int_{\mathcal{T}_k} p(\mathbf{t}_k \mid \theta) \log[p(\mathbf{t}_k \mid \theta)] \, d\mathbf{t}_k \, d\theta$$

$$- \int_{\mathcal{T}_k} p_0(\mathbf{t}_k) \log[p_0(\mathbf{t}_k)] \, d\mathbf{t}_k.$$

Because $I\{p_0 \mid \mathcal{M}^k\} \geq 0$,

$$\int_{\mathcal{T}_k} p_0(\mathbf{t}_k) \log[p_0(\mathbf{t}_k)] \, d\mathbf{t}_k \leq \int_{\Theta_0} p_0(\theta) \int_{\mathcal{T}_k} p(\mathbf{t}_k \mid \theta) \log[p(\mathbf{t}_k \mid \theta)] \, d\mathbf{t}_k \, d\theta.$$

Condition (i) and the continuity of $p$ ensure the right-hand side of the last equation is bounded above, and condition (ii) ensures that its left-hand side is bounded below. Consequently, $I\{p_0 \mid \mathcal{M}^k\} < \infty$.

## APPENDIX F: PROOF OF THEOREM 7

For any $p(\theta) \in \mathcal{P}_s$, denote the posterior corresponding to $p_0$ (the restriction of $p$ to the compact set $\Theta_0$) by $p_0(\theta \mid \mathbf{t}_k)$.

*Step* 1. We give an expansion of $I\{p_0 \mid \mathcal{M}^k\}$, defined by

$$I\{p_0 \mid \mathcal{M}^k\} = \int_{\Theta_0} p_0(\theta) \int_{\mathcal{T}_k} p(\mathbf{t}_k \mid \theta) \log\left[\frac{p_0(\theta \mid \mathbf{t}_k)}{p_0(\theta)}\right] d\mathbf{t}_k \, d\theta.$$

Use the equality

$$\frac{p_0(\theta \mid \mathbf{t}_k)}{p_0(\theta)} = \frac{p_0(\theta \mid \mathbf{t}_k)}{\pi_0^*(\theta \mid \mathbf{t}_k)} \frac{\pi_0^*(\theta \mid \mathbf{t}_k)}{\pi^*(\theta \mid \mathbf{t}_k)} \frac{\pi^*(\theta \mid \mathbf{t}_k)}{\pi_{0k}^*(\theta)} \frac{\pi_{0k}^*(\theta)}{p_0(\theta)},$$

where

(F.1) $\quad \pi_{0k}^*(\theta) = \dfrac{f_k(\theta)}{c_0(f_k)} 1_{\Theta_0}(\theta) \quad \text{and} \quad c_0(f_k) = \int_{\Theta_0} f_k(\theta) \, d\theta.$

We have the decomposition

(F.2) $$I\{p_0 \mid \mathcal{M}^k\} = \sum_{j=1}^{4} G_{jk},$$

where

$$G_{1k} = -\int_{\Theta_0} p_0(\theta) \int_{\mathcal{T}_k} p(\mathbf{t}_k \mid \theta) \log\left[\frac{\pi_0^*(\theta \mid \mathbf{t}_k)}{p_0(\theta \mid \mathbf{t}_k)}\right] d\theta \, d\mathbf{t}_k,$$

$$G_{2k} = \int_{\Theta_0} p_0(\theta) \int_{\mathcal{T}_k} p(\mathbf{t}_k \mid \theta) \log\left[\frac{\pi_0^*(\theta \mid \mathbf{t}_k)}{\pi^*(\theta \mid \mathbf{t}_k)}\right] d\mathbf{t}_k \, d\theta,$$

$$G_{3k} = \int_{\Theta_0} p_0(\theta) \int_{\mathcal{T}_k} p(\mathbf{t}_k \mid \theta) \log\left[\frac{\pi^*(\theta \mid \mathbf{t}_k)}{\pi_{0k}^*(\theta)}\right] d\mathbf{t}_k \, d\theta,$$

$$G_{4k} = \int_{\Theta_0} p_0(\theta) \int_{\mathcal{T}_k} p(\mathbf{t}_k \mid \theta) \log\left[\frac{\pi_{0k}^*(\theta)}{p_0(\theta)}\right] d\mathbf{t}_k \, d\theta.$$

It is easy to see that

$$G_{3k} = \int_{\Theta_0} p_0(\theta) \log\left[\frac{f_k(\theta)}{\pi_{0k}^*(\theta)}\right] d\theta.$$

From (F.1), $f_k(\theta)/\pi_{0k}^*(\theta) = c_0(f_k)$ on $\Theta_0$. Then,

(F.3) $$G_{3k} = \log[c_0(f_k)].$$

Clearly,

(F.4) $$G_{4k} = -\int_{\Theta_0} p_0(\theta) \log\left[\frac{p_0(\theta)}{\pi_{0k}^*(\theta)}\right] d\theta.$$

Note that the continuity of $p(\mathbf{t}_k \mid \theta)$ in $\theta$ and integrability will imply the continuity of $f_k$. So, $\pi_{0k}^*$ is continuous and bounded, and $G_{4k}$ is finite. Since, $0 \leq I\{p_0 \mid \mathcal{M}^k\} < \infty$, $G_{jk}$, $j = 1, 2, 3$ are all nonnegative and finite.

*Step* 2. We show that

(F.5) $$\lim_{k \to \infty} G_{1k} = 0 \qquad \forall p \in \mathcal{P}_s.$$

It is easy to see

$$\frac{\pi_0^*(\theta \mid \mathbf{t}_k)}{p_0(\theta \mid \mathbf{t}_k)} = \frac{\pi^*(\theta)}{p(\theta)} \frac{\int_{\Theta_0} p(\mathbf{t}_k \mid \tau) p(\tau) \, d\tau}{\int_{\Theta_0} p(\mathbf{t}_k \mid \tau) \pi^*(\tau) \, d\tau}$$

(F.6) $$= \frac{\pi^*(\theta)}{p(\theta)} \frac{\int_{\Theta_0} p(\mathbf{t}_k \mid \tau) p(\tau) \, d\tau}{\int_{\Theta} p(\mathbf{t}_k \mid \tau) \pi^*(\tau) \, d\tau} [P^*(\Theta_0 \mid \mathbf{t}_k)]^{-1}.$$

The definition of posterior consistency of $\pi^*$ is that, for any $\theta \in \Theta$ and any $\varepsilon > 0$,

(F.7) $$P^*(|\tau - \theta| \leq \varepsilon \mid \mathbf{t}_k) \equiv \int_{\{\tau : |\tau - \theta| \leq \varepsilon\}} \pi^*(\tau \mid \mathbf{t}_k) \, d\tau \xrightarrow{P} 1,$$

in probability $p(\mathbf{t}_k \mid \theta)$ as $k \to \infty$. It is immediate that

(F.8) $$P^*(\Theta_0 \mid \mathbf{t}_k) = \frac{\int_{\Theta_0} p(\mathbf{t}_k \mid \tau) \pi^*(\tau) \, d\tau}{\int_{\Theta} p(\mathbf{t}_k \mid \tau) \pi^*(\tau) \, d\tau} \xrightarrow{P} 1,$$

with probability $p(\mathbf{t}_k \mid \theta)$ as $\theta \in \Theta_0$ and $k \to \infty$. Because both $\pi^*$ and $p$ are continuous, for any $\varepsilon > 0$, there is small $\delta > 0$, such that

$$(\text{F.9}) \qquad \left| \frac{p(\tau)}{\pi^*(\tau)} - \frac{p(\theta)}{\pi^*(\theta)} \right| \leq \varepsilon \qquad \forall \tau \in \Theta_0 \cap (\theta - \delta, \theta + \delta)$$

For such a $\delta$, we could write

$$(\text{F.10}) \qquad \frac{\pi^*(\theta)}{p(\theta)} \frac{\int_{\Theta_0} p(\mathbf{t}_k \mid \tau) p(\tau) \, d\tau}{\int_\Theta p(\mathbf{t}_k \mid \tau) \pi^*(\tau) \, d\tau} \equiv J_{1k} + J_{2k},$$

where

$$J_{1k} = \frac{\pi^*(\theta)}{p(\theta)} \frac{\int_{\Theta_0 \cap (\theta-\delta,\theta+\delta)} p(\mathbf{t}_k \mid \tau)(p(\tau)/\pi^*(\tau))\pi^*(\tau) \, d\tau}{\int_\Theta p(\mathbf{t}_k \mid \tau) \pi^*(\tau) \, d\tau},$$

$$J_{2k} = \frac{\pi^*(\theta)}{p(\theta)} \frac{\int_{\Theta_0 \cap (\theta-\delta,\theta+\delta)^c} p(\mathbf{t}_k \mid \tau)(p(\tau)/\pi^*(\tau))\pi^*(\tau) \, d\tau}{\int_\Theta p(\mathbf{t}_k \mid \tau) \pi^*(\tau) \, d\tau}.$$

Clearly, (F.9) implies that

$$J_{1k} \geq \frac{\pi^*(\theta)}{p(\theta)} \left[ \frac{p(\theta)}{\pi^*(\theta)} - \varepsilon \right] \int_{\Theta_0 \cap (\theta-\delta,\theta+\delta)} \pi^*(\tau \mid \mathbf{t}_k) \, d\tau,$$

$$J_{1k} \leq \frac{\pi^*(\theta)}{p(\theta)} \left[ \frac{p(\theta)}{\pi^*(\theta)} + \varepsilon \right] \int_{\Theta_0 \cap (\theta-\delta,\theta+\delta)} \pi^*(\tau \mid \mathbf{t}_k) \, d\tau.$$

(F.7) implies that, for the fixed $\delta$ and $\theta \in \Theta_0$,

$$(\text{F.11}) \qquad \left[ 1 - \varepsilon \frac{\pi^*(\theta)}{p(\theta)} \right] \leq J_{1k} \leq \left[ 1 + \varepsilon \frac{\pi^*(\theta)}{p(\theta)} \right]$$

with probability $p(\mathbf{t}_k \mid \theta)$ as $k \to \infty$. Noting that $p(\theta)$ is continuous and positive on $\Theta_0$, let $M_1 > 0$ and $M_2$ be the lower and upper bounds of $p$ on $\Theta_0$. From (F.7),

$$(\text{F.12}) \qquad 0 \leq J_{2k} \leq \frac{M_2 \pi^*(\theta)}{M_1 p(\theta)} \int_{\Theta_0 \cap (\theta-\delta,\theta+\delta)^c} \pi^*(\tau \mid \mathbf{t}_k) \, d\tau \xrightarrow{P} 0,$$

with probability $p(\mathbf{t}_k \mid \theta)$ as $k \to \infty$. Combining (F.6), (F.8) and (F.10)–(F.12), we know that

$$(\text{F.13}) \qquad \frac{\pi_0^*(\theta \mid \mathbf{t}_k)}{p_0(\theta \mid \mathbf{t}_k)} \xrightarrow{P} 1$$

with probability $p(\mathbf{t}_k \mid \theta)$ as $k \to \infty$. It is easy to see that the left quantity of (F.13) is bounded above and below, so the dominated convergence theorem implies (F.5).

*Step* 3. We show that

$$(\text{F.14}) \quad G_{5k} \equiv \int_{\Theta_0} \pi_0^*(\theta) \int_{\mathcal{T}_k} p(\mathbf{t}_k \mid \theta) \log \frac{\pi_0^*(\theta \mid \mathbf{t}_k)}{\pi^*(\theta \mid \mathbf{t}_k)} \, d\mathbf{t}_k \, d\theta \to 0 \qquad \text{as } k \to \infty.$$

For any measurable set $A \subset \mathbb{R}$, denote $P^*(A \mid \mathbf{t}_k) = \int_A \pi^*(\theta \mid \mathbf{t}_k)\, d\theta$. Then,

$$\frac{\pi_0^*(\theta \mid \mathbf{t}_k)}{\pi^*(\theta \mid \mathbf{t}_k)} = \frac{p(\mathbf{t}_k \mid \theta)\pi_0^*(\theta)/p_0^*(\mathbf{t}_k)}{p(\mathbf{t}_k \mid \theta)\pi^*(\theta)/p^*(\mathbf{t}_k)} = \frac{\int_\Theta p(\mathbf{t}_k \mid \theta)\pi^*(\theta)\, d\theta}{\int_{\Theta_0} p(\mathbf{t}_k \mid \theta)\pi^*(\theta)\, d\theta}$$

$$= 1 + \frac{\int_{\Theta_0^c} p(\mathbf{t}_k \mid \theta)\pi^*(\theta)\, d\theta}{\int_{\Theta_0} p(\mathbf{t}_k \mid \theta)\pi^*(\theta)\, d\theta} = 1 + \frac{P^*(\Theta_0^c \mid \mathbf{t}_k)}{P^*(\Theta_0 \mid \mathbf{t}_k)}.$$

Thus,

$$G_{5k} = \int_{\Theta_0} \pi_0^*(\theta) \int_{\mathcal{T}_k} p(\mathbf{t}_k \mid \theta) \log\left\{1 + \frac{P^*(\Theta_0^c \mid \mathbf{t}_k)}{P^*(\Theta_0 \mid \mathbf{t}_k)}\right\} d\mathbf{t}_k\, d\theta.$$

For any $0 \leq a \leq b \leq \infty$, denote

(F.15) $$\mathcal{T}_{k,a,b} = \left\{\mathbf{t}_k : a \leq \frac{P^*(\Theta_0^c \mid \mathbf{t}_k)}{P^*(\Theta_0 \mid \mathbf{t}_k)} < b\right\}.$$

Clearly, if $0 < \varepsilon < M < \infty$,

$$\mathcal{T}_k = \mathcal{T}_{k,0,\varepsilon} \cup \mathcal{T}_{k,\varepsilon,M} \cup \mathcal{T}_{k,M,\infty}.$$

We then have the decomposition for $G_{5k}$,

(F.16) $$G_{5k} \equiv G_{5k1} + G_{5k2} + G_{5k3},$$

where

$$G_{5k1} = \int_{\Theta_0} \pi_0^*(\theta) \int_{\mathcal{T}_{k,0,\varepsilon}} p(\mathbf{t}_k \mid \theta) \log\left\{1 + \frac{P^*(\Theta_0^c \mid \mathbf{t}_k)}{P^*(\Theta_0 \mid \mathbf{t}_k)}\right\} d\mathbf{t}_k\, d\theta,$$

$$G_{5k2} = \int_{\Theta_0} \pi_0^*(\theta) \int_{\mathcal{T}_{k,\varepsilon,M}} p(\mathbf{t}_k \mid \theta) \log\left\{1 + \frac{P^*(\Theta_0^c \mid \mathbf{t}_k)}{P^*(\Theta_0 \mid \mathbf{t}_k)}\right\} d\mathbf{t}_k\, d\theta,$$

$$G_{5k3} = \int_{\Theta_0} \pi_0^*(\theta) \int_{\mathcal{T}_{k,M,\infty}} p(\mathbf{t}_k \mid \theta) \log\left\{1 + \frac{P^*(\Theta_0^c \mid \mathbf{t}_k)}{P^*(\Theta_0 \mid \mathbf{t}_k)}\right\} d\mathbf{t}_k\, d\theta.$$

The posterior consistency (F.8) implies that if $\theta \in \Theta_0^0$ (the interior of $\Theta_0$),

(F.17) $$P^*(\Theta_0^c \mid \mathbf{t}_k) \xrightarrow{P} 0,$$

in probability $p(\mathbf{t}_k \mid \theta)$ as $k \to \infty$. So (F.17) implies that, for any small $\varepsilon > 0$ and any fixed $\theta \in \Theta_0^0$,

(F.18) $$\int_{\mathcal{T}_{k,\varepsilon,\infty}} p(\mathbf{t}_k \mid \theta)\, d\mathbf{t}_k \longrightarrow 0 \quad \text{as } k \to \infty.$$

For small $\varepsilon > 0$,

$$G_{5k1} \leq \log(1 + \varepsilon) \int_{\Theta_0} \pi_0^*(\theta) \int_{\mathcal{T}_{k,0,\varepsilon}} p(\mathbf{t}_k \mid \theta)\, d\mathbf{t}_k\, d\theta$$

$$\leq \log(1 + \varepsilon) < \varepsilon.$$

For any large $M > \max(\varepsilon, e - 1)$,

$$G_{5k2} \leq \log(1 + M) \int_{\Theta_0} \pi_0^*(\theta) \int_{\mathcal{T}_{k,\varepsilon,M}} p(\mathbf{t}_k \mid \theta) \, d\mathbf{t}_k \, d\theta$$

$$\leq \log(1 + M) \int_{\Theta_0} \pi_0^*(\theta) \int_{\mathcal{T}_{k,\varepsilon,\infty}} p(\mathbf{t}_k \mid \theta) \, d\mathbf{t}_k \, d\theta.$$

Since $\pi_0^*$ is bounded on $\Theta_0$, (F.18) and dominated convergence theorem imply that

$$G_{5k2} \to 0 \qquad \text{as } k \to \infty.$$

Also,

$$G_{5k3} = \int_{\Theta_0} \pi_0^*(\theta) \int_{\mathcal{T}_{k,M,\infty}} p(\mathbf{t}_k \mid \theta) \log\left\{\frac{1}{P^*(\Theta_0 \mid \mathbf{t}_k)}\right\} d\mathbf{t}_k \, d\theta$$

$$= -\frac{1}{c_0(\pi^*)} \int_{\mathcal{T}_{k,M,\infty}} p^*(\mathbf{t}_k) \int_{\Theta_0} p^*(\theta \mid \mathbf{t}_k) \log[P^*(\Theta_0 \mid \mathbf{t}_k)] \, d\theta \, d\mathbf{t}_k$$

$$= -\frac{1}{c_0(\pi^*)} \int_{\mathcal{T}_{k,M,\infty}} p^*(\mathbf{t}_k) P^*(\Theta_0 \mid \mathbf{t}_k) \log[P^*(\Theta_0 \mid \mathbf{t}_k)] \, d\mathbf{t}_k.$$

Note that $\mathbf{t}_k \in \mathcal{T}_{k,M,\infty}$ if and only if $P^*(\Theta_0 \mid \mathbf{t}_k) < 1/(1+M)$. Also, $-p \log(p)$ is increasing for $p \in (0, 1/e)$. This implies that

$$-P^*(\Theta_0 \mid \mathbf{t}_k) \log[P^*(\Theta_0 \mid \mathbf{t}_k)] < \frac{1}{1 + M} \log(1 + M).$$

Consequently,

$$J_{5k} \leq \frac{1}{c_0(\pi^*)(1 + M)} \log(1 + M) \int_{\mathcal{T}_{k,M,\infty}} p^*(\mathbf{t}_k) \, d\mathbf{t}_k$$

$$\leq \frac{1}{c_0(\pi^*)(1 + M)} \log(1 + M).$$

Now for fixed small $\varepsilon > 0$, we could choose $M > \max(\varepsilon, e - 1)$ large enough so that $G_{5k3} \leq \varepsilon$. For such fixed $\varepsilon$ and $M$, we know $G_{5k2} \to 0$ as $k \to \infty$. Since $\varepsilon$ is arbitrary, (F.14) holds.

*Step* 4. We show that

(F.19) $$\lim_{k \to \infty} G_{2k} = 0 \qquad \forall p \in \mathcal{P}_s.$$

Note that for any $p \in \mathcal{P}_s$, there is a constant $M > 0$, such that

$$\sup_{\tau \in \Theta_0} \frac{p_0(\tau)}{\pi_0^*(\tau)} \leq M.$$

Since $\pi_0^*(\theta \mid \mathbf{t}_k)/\pi^*(\theta \mid \mathbf{t}_k) \geq 1$,

$$0 \leq G_{2k} \leq M \int_{\Theta_0} \pi_0^*(\theta) \int_{\mathcal{T}_k} p(\mathbf{t}_k \mid \theta) \log\left[\frac{\pi_0^*(\theta \mid \mathbf{t}_k)}{\pi^*(\theta \mid \mathbf{t}_k)}\right] d\mathbf{t}_k \, d\theta = MG_{5k}.$$

Then, (F.14) implies (F.19) immediately.

*Step* 5. It follows from (F.2) that for any prior $p \in \mathcal{P}_s$,

$$I\{\pi_0 \mid \mathcal{M}^k\} - I\{p_0 \mid \mathcal{M}^k\}$$

$$= -G_{1k} - G_{2k} + \int_{\Theta_0} \pi_0(\theta) \int_{\mathcal{T}_k} p(\mathbf{t}_k \mid \theta) \log\left[\frac{\pi_0^*(\theta \mid \mathbf{t}_k)}{\pi^*(\theta \mid \mathbf{t}_k)}\right] d\mathbf{t}_k \, d\theta$$

$$- \int_{\Theta_0} \pi_0(\theta) \log\left[\frac{\pi_0(\theta)}{\pi_{0k}^*(\theta)}\right] d\theta + \int_{\Theta_0} p_0(\theta) \log\left[\frac{p_0(\theta)}{\pi_{0k}^*(\theta)}\right] d\theta.$$

Steps 2 and 4 imply that

$$\lim_{k \to \infty} (I\{\pi_0 \mid \mathcal{M}^k\} - I\{p_0 \mid \mathcal{M}^k\}) = \lim_{k \to \infty} \left\{ -\int_{\Theta_0} \pi_0(\theta) \log\left[\frac{\pi_0(\theta)}{\pi_{0k}^*(\theta)}\right] d\theta \right.$$

(F.20) $$\left. + \int_{\Theta_0} p_0(\theta) \log\left[\frac{p_0(\theta)}{\pi_{0k}^*(\theta)}\right] d\theta \right\}$$

$$\geq -\lim_{k \to \infty} \int_{\Theta_0} \pi_0(\theta) \log\left[\frac{\pi_0(\theta)}{\pi_{0k}^*(\theta)}\right] d\theta,$$

the last inequality holding since the second term is always nonnegative. Finally,

$$\lim_{k \to \infty} \int_{\Theta_0} \pi_0(\theta) \log[\pi_{0k}^*(\theta)] \, d\theta$$

$$= \lim_{k \to \infty} \int_{\Theta_0} \pi_0(\theta) \log\left[\frac{f_k(\theta)}{c_0(f_k)}\right] d\theta$$

$$= \lim_{k \to \infty} \int_{\Theta_0} \pi_0(\theta) \log\left[\frac{f_k(\theta)}{f_k(\theta_0)} \frac{f_k(\theta_0)}{c_0(f_k)}\right] d\theta$$

$$= \lim_{k \to \infty} \int_{\Theta_0} \pi_0(\theta) \log\left[\frac{f_k(\theta)}{f_k(\theta_0)}\right] d\theta + \lim_{k \to \infty} \log\left[\frac{f_k(\theta_0)}{c_0(f_k)}\right]$$

$$= \int_{\Theta_0} \pi_0(\theta) \log[f(\theta)] \, d\theta - \log[c_0(f)]$$

$$= \int_{\Theta_0} \pi_0(\theta) \log[\pi_0(\theta)] \, d\theta,$$

the second to last line following from condition (i) and

$$\lim_{k \to \infty} \frac{c_0(f_k)}{f_k(\theta_0)} = \lim_{k \to \infty} \int_{\Theta_0} \frac{f_k(\theta)}{f_k(\theta_0)} \, d\theta$$

$$= \int_{\Theta_0} \lim_{k \to \infty} \frac{f_k(\theta)}{f_k(\theta_0)} d\theta = \int_{\Theta_0} f(\theta) d\theta = c_0(f).$$

Consequently, the right-hand side of (F.20) is 0, completing the proof.

## APPENDIX G: PROOF OF THEOREM 8

Let $\mathbf{x}^{(k)} = \{x_1, \ldots, x_k\}$ consist of $k$ replications from the original uniform distribution on the interval $(a_1(\theta), a_2(\theta))$. Let $t_1 = t_{k1} = \min\{x_1, \ldots, x_k\}$ and $t_2 = t_{k2} = \max\{x_1, \ldots, x_k\}$. Clearly, $\mathbf{t}_k \equiv (t_1, t_2)$ are sufficient statistics with density

$$(G.1) \quad p(t_1, t_2 \mid \theta) = \frac{k(k-1)(t_2 - t_1)^{k-2}}{[a_2(\theta) - a_1(\theta)]^k}, \quad a_1(\theta) < t_1 < t_2 < a_2(\theta).$$

Choosing $\pi^*(\theta) = 1$, the corresponding posterior density of $\theta$ is

$$(G.2) \quad \pi^*(\theta \mid t_1, t_2) = \frac{1}{[a_2(\theta) - a_1(\theta)]^k m_k(t_1, t_2)},$$

$$a_2^{-1}(t_2) < \theta < a_1^{-1}(t_1),$$

where

$$(G.3) \quad m_k(t_1, t_2) = \int_{a_2^{-1}(t_2)}^{a_1^{-1}(t_1)} \frac{1}{[a_2(s) - a_1(s)]^k} ds.$$

Consider the transformation

$$(G.4) \quad y_1 = k(a_1^{-1}(t_1) - \theta) \quad \text{and} \quad y_2 = k(\theta - a_2^{-1}(t_2)),$$

or equivalently, $t_1 = a_1(\theta + y_1/k)$ and $t_2 = a_2(\theta - y_2/k)$.

We first consider the frequentist asymptotic distribution of $(y_1, y_2)$. For $\theta > \theta_0$, we know $a_1(\theta) < a_2(\theta)$. For any fixed $y_1 > 0$ and $y_2 > 0$, $a_1(\theta + y_1/k) < a_2(\theta - y_2/k)$ when $k$ is large enough. From (G.1), the joint density of $(y_1, y_2)$ is

$$p(y_1, y_2 \mid \theta)$$
$$= \frac{(k-1)}{k} \frac{a_1'(\theta + y_1/k) a_2'(\theta - y_2/k)}{[a_2(\theta) - a_1(\theta)]^k} \left\{ a_2\left(\theta - \frac{y_2}{k}\right) - a_1\left(\theta + \frac{y_1}{k}\right) \right\}^{k-2}$$
$$= \frac{(k-1)}{k} \frac{a_1'(\theta + y_1/k) a_2'(\theta - y_2/k)}{[a_2(\theta) - a_1(\theta)]^2}$$
$$\times \left\{ 1 - \frac{a_1'(\theta) y_1 + a_2'(\theta) y_2}{[a_2(\theta) - a_1(\theta)]k} + o\left(\frac{1}{k}\right) \right\}^{k-2}.$$

For fixed $\theta > \theta_0$, $y_1, y_2 > 0$, as $k \to \infty$,

$$p(y_1, y_2 \mid \theta) \to \frac{a_1'(\theta) a_2'(\theta)}{[a_2(\theta) - a_1(\theta)]^2} \exp\left\{-\frac{a_1'(\theta) y_1 + a_2'(\theta) y_2}{a_2(\theta) - a_1(\theta)}\right\}$$

(G.5)
$$\equiv p^*(y_1, y_2 \mid \theta).$$

Consequently, as $k \to \infty$, the $y_i$'s have independent exponential distributions with means $\lambda_i = [a_2(\theta) - a_1(\theta)]/a_i'(\theta)$.

With the transformation (G.4),

$$m_k(t_1, t_2) = \int_{\theta - y_2/k}^{\theta + y_1/k} \frac{1}{[a_2(s) - a_1(s)]^k} \, ds$$

(G.6)
$$= \frac{1}{k} \int_{-y_2}^{y_1} \frac{1}{[a_2(\theta + v/k) - a_1(\theta + v/k)]^k} \, dv.$$

So, for any fixed $y_1, y_2 > 0$ as $k \to \infty$,

$$k[a_2(\theta) - a_1(\theta)]^k m_k(t_1, t_2)$$

$$\longrightarrow \int_{-y_2}^{y_1} \exp\left[-\frac{a_2'(\theta) - a_1'(\theta)}{a_2(\theta) - a_1(\theta)} v\right] dv$$

$$= \frac{a_2(\theta) - a_1(\theta)}{a_2'(\theta) - a_1'(\theta)} \exp\left[\frac{a_2'(\theta) - a_1'(\theta)}{a_2(\theta) - a_1(\theta)} y_2\right]$$

$$\times \left\{1 - \exp\left[-\frac{a_2'(\theta) - a_1'(\theta)}{a_2(\theta) - a_1(\theta)} (y_1 + y_2)\right]\right\}.$$

Then, for fixed $\theta > \theta_0$ as $k \to \infty$,

$$\int \log(\pi(\theta \mid t_1, t_2)) f(t_1, t_2 \mid \theta) \, dt_1 \, dt_2 - \log(k)$$

(G.7)
$$\longrightarrow \log\left\{\frac{a_2'(\theta) - a_1'(\theta)}{a_2(\theta) - a_1(\theta)}\right\} + J_1(\theta) + J_2(\theta),$$

where

$$J_1(\theta) = \frac{a_2'(\theta) - a_1'(\theta)}{a_2(\theta) - a_1(\theta)} \int_0^\infty \int_0^\infty y_2 p^*(y_1, y_2 \mid \theta) \, dy_1 \, dy_2,$$

$$J_2(\theta) = -\int_0^\infty \int_0^\infty \log\left\{1 - \exp\left[-\frac{a_2'(\theta) - a_1'(\theta)}{a_2(\theta) - a_1(\theta)} (y_1 + y_2)\right]\right\}$$
$$\times p^*(y_1, y_2 \mid \theta) \, dy_1 \, dy_2.$$

It follows from (G.5) that

$$J_1(\theta) = -b_2,$$
$$J_2(\theta) = -E \log\{1 - e^{-b_1 V_1} e^{-b_2 V_2}\},$$

DEFINITION OF REFERENCE PRIORS 33

where $V_1$ and $V_2$ are i.i.d. with the standard exponential distribution. Then,

$$J_2(\theta) = \sum_{j=1}^{\infty} \frac{1}{j} E(e^{-jb_1 V_1}) E(e^{-jb_2 V_2})$$

(G.8)
$$= \sum_{j=1}^{\infty} \frac{1}{j(b_1 j + 1)(b_2 j + 1)}$$

$$= \frac{1}{b_1 - b_2} \sum_{j=1}^{\infty} \frac{1}{j} \left( \frac{1}{j + 1/b_1} - \frac{1}{j + 1/b_2} \right).$$

Note that the digamma function $\psi(z)$ satisfies the equation,

(G.9)
$$\sum_{j=1}^{\infty} \frac{1}{j(j+z)} = \frac{\psi(z+1) + \gamma}{z},$$

for $z > 0$, where $\gamma$ is the Euler–Mascherono constant (see, e.g., Boros and Moll [14].) Equations (G.8) and (G.9) imply that

$$J_2(\theta) = \frac{1}{b_1 - b_2} \left\{ b_1 \left[ \psi\left(\frac{1}{b_1} + 1\right) + \gamma \right] - b_2 \left[ \psi\left(\frac{1}{b_2} + 1\right) + \gamma \right] \right\}$$

$$= \gamma + \frac{1}{b_1 - b_2} \left\{ b_1 \psi\left(\frac{1}{b_1} + 1\right) - b_2 \psi\left(\frac{1}{b_2} + 1\right) \right\}.$$

Using the fact that $\psi(z+1) = \psi(z) + 1/z$,

(G.10)
$$J_1(\theta) + J_2(\theta) = \gamma - b_2 + \frac{1}{b_1 - b_2} \left\{ b_1^2 + b_1 \psi\left(\frac{1}{b_1}\right) - b_2^2 - b_2 \psi\left(\frac{1}{b_2}\right) \right\}$$

$$= \gamma + b_1 + \frac{1}{b_1 - b_2} \left\{ b_1 \psi\left(\frac{1}{b_1}\right) - b_2 \psi\left(\frac{1}{b_2}\right) \right\}.$$

The result follows from (G.7) and (G.10).

**Acknowledgments.** The authors are grateful to Susie Bayarri for helpful discussions. The authors acknowledge the constructive comments of the Associate Editor and referee.

J. O. BERGER
DEPARTMENT OF STATISTICAL SCIENCES
DUKE UNIVERSITY
BOX 90251
DURHAM, NORTH CAROLINA 27708-0251
USA
E-MAIL: berger@stat.duke.edu
URL: www.stat.duke.edu/~berger/

J. M. BERNARDO
DEPARTAMENTO DE ESTADÍSTICA
FACULTAD DE MATEMÁTICAS
46100–BURJASSOT
VALENCIA
SPAIN
E-MAIL: jose.m.bernardo@uv.es
URL: www.uv.es/~bernardo/

D. Sun
Department of Statistics
University of Missouri-Columbia
146 Middlebush Hall
Columbia, Missouri 65211-6100
USA
E-mail: sund@missouri.edu
URL: www.stat.missouri.edu/~dsun/